\documentclass[12pt]{article}

\usepackage{fontspec, xunicode, xltxtra}  
\usepackage{mathtools}
\usepackage{amsmath,diagbox}
\usepackage{indentfirst}
\usepackage{mathrsfs}
\usepackage{amsfonts}
\usepackage{arydshln}
\usepackage{enumerate}
\usepackage{setspace}
\usepackage{amssymb,amsthm,cases}
\usepackage{amsbsy,pifont}
\usepackage{latexsym,diagbox,tabularx}
\usepackage{amsmath,amscd,bbm,multirow}
\usepackage{graphicx,epsfig,extarrows,dsfont,mathtools}
\usepackage[final,allcolors=blue,colorlinks=true]{hyperref}
\usepackage[normalem]{ulem}
\input xy
\usepackage{caption,subcaption} 
\usepackage{tikz}
\usepackage{tikz-cd}
\usetikzlibrary{decorations.pathreplacing,decorations.markings}
\usetikzlibrary{cd}
\usetikzlibrary{shapes.geometric}
\usetikzlibrary{arrows,arrows.meta}
\usetikzlibrary{graphs,positioning}
\usetikzlibrary{matrix,decorations.pathmorphing}
\usetikzlibrary{math,calc,intersections,through,angles,arrows.meta,shapes.geometric,shadows,quotes,spy,datavisualization,datavisualization.formats.functions,plotmarks}
\tikzset{every picture/.style={samples=300,smooth,line join=round,thick,>=stealth}}
\usepackage{enumerate}
\usepackage{fancyhdr}
\usepackage[sectionbib]{chapterbib}

\tikzset{
	on each segment/.style={
		decorate,
		decoration={
			show path construction,
			moveto code={},
			lineto code={
				\path [#1]
				(\tikzinputsegmentfirst) -- (\tikzinputsegmentlast);
			},
			curveto code={
				\path [#1] (\tikzinputsegmentfirst)
				.. controls
				(\tikzinputsegmentsupporta) and (\tikzinputsegmentsupportb)
				..
				(\tikzinputsegmentlast);
			},
			closepath code={
				\path [#1]
				(\tikzinputsegmentfirst) -- (\tikzinputsegmentlast);
			},
		},
	},
	mid arrow/.style={postaction={decorate,decoration={
				markings,
				mark=at position .5 with {\arrow[#1]{stealth}}
	}}},
}

\numberwithin{equation}{section}

\theoremstyle{plain}
\newtheorem{theorem}{Theorem}[section]

\newtheorem{corollary}[theorem]{Corollary}

\newtheorem{lemma}[theorem]{Lemma}

\theoremstyle{definition}
\newtheorem{definition}[theorem]{Definition}

\newtheorem{remark}[theorem]{Remark}

\def\XXint#1#2#3{{\setbox0=\hbox{$#1{#2#3}{\int}$}
		\vcenter{\hbox{$#2#3$}}\kern-.5\wd0}}

\DeclareMathSymbol{\subseteq}{\mathrel}{symbols}{"12}
\DeclareMathSymbol{\supseteq}{\mathrel}{symbols}{"13} 
\DeclareMathSymbol{\subsetneq}{\mathrel}{AMSb}{"28}                  \DeclareMathSymbol{\supsetneq}{\mathrel}{AMSb}{"29}    

\DeclareMathSymbol{\nsubseteq}{\mathrel}{AMSb}{"2A}                  

\DeclareMathSymbol{\nsupseteq}{\mathrel}{AMSb}{"2B}

\DeclareMathDelimiter{\langle}{\mathop}{symbols}{"68}{largesymbols}{"0A}
\DeclareMathDelimiter{\rangle}{\mathclose}{symbols}{"69}{largesymbols}{"0B}

\DeclareSymbolFont{txfontsA}{U}{txmia}{m}{it}
\SetSymbolFont{txfontsA}{bold}{U}{txmia}{bx}{it}
\DeclareFontSubstitution{U}{txmia}{m}{it}
\DeclareMathSymbol{\upalpha}{\mathord}{txfontsA}{"0B}
\DeclareMathSymbol{\upbeta}{\mathord}{txfontsA}{"0C}
\DeclareMathSymbol{\upgamma}{\mathord}{txfontsA}{"0D}
\DeclareMathSymbol{\updelta}{\mathord}{txfontsA}{"0E}
\DeclareMathSymbol{\upepsilon}{\mathord}{txfontsA}{"0F}
\DeclareMathSymbol{\upzeta}{\mathord}{txfontsA}{"10}
\DeclareMathSymbol{\upeta}{\mathord}{txfontsA}{"11}
\DeclareMathSymbol{\uptheta}{\mathord}{txfontsA}{"12}
\DeclareMathSymbol{\upiota}{\mathord}{txfontsA}{"13}
\DeclareMathSymbol{\upkappa}{\mathord}{txfontsA}{"14}
\DeclareMathSymbol{\uplambda}{\mathord}{txfontsA}{"15}
\DeclareMathSymbol{\upmu}{\mathord}{txfontsA}{"16}
\DeclareMathSymbol{\upnu}{\mathord}{txfontsA}{"17}
\DeclareMathSymbol{\upxi}{\mathord}{txfontsA}{"18}
\DeclareMathSymbol{\uppi}{\mathord}{txfontsA}{"19}
\DeclareMathSymbol{\uprho}{\mathord}{txfontsA}{"1A}
\DeclareMathSymbol{\upsigma}{\mathord}{txfontsA}{"1B}
\DeclareMathSymbol{\uptau}{\mathord}{txfontsA}{"1C}
\DeclareMathSymbol{\upupsilon}{\mathord}{txfontsA}{"1D}
\DeclareMathSymbol{\upphi}{\mathord}{txfontsA}{"1E}
\DeclareMathSymbol{\upchi}{\mathord}{txfontsA}{"1F}
\DeclareMathSymbol{\uppsi}{\mathord}{txfontsA}{"20}
\DeclareMathSymbol{\upomega}{\mathord}{txfontsA}{"21}
\DeclareMathSymbol{\upvarepsilon}{\mathord}{txfontsA}{"22}
\DeclareMathSymbol{\upvartheta}{\mathord}{txfontsA}{"23}
\DeclareMathSymbol{\upvarpi}{\mathord}{txfontsA}{"24}
\DeclareMathSymbol{\upvarrho}{\mathord}{txfontsA}{"25}
\DeclareMathSymbol{\upvarsigma}{\mathord}{txfontsA}{"26}
\DeclareMathSymbol{\upvarphi}{\mathord}{txfontsA}{"27}                   

\DeclareSymbolFont{ugmL}{OMX}{mdugm}{m}{n}
\SetSymbolFont{ugmL}{bold}{OMX}{mdugm}{b}{n}
\DeclareMathAccent{\wideparen}{\mathord}{ugmL}{"F3}

\def\pt{\partial}

\def\ra{\rightarrow}

\def\s{\subseteq}

\def\e{\epsilon}

\def\ol{\overline}

\def\vp{\varphi}
\def\lg{\langle}

\def\rg{\rangle}

\def\bf{\textbf}
\def\pt{\partial}

\def\Om{\Omega}
\def\la{\lambda}
\def\al{\alpha}
\def\be{\beta}
\def\de{\delta}
\def\ga{\gamma}

\def\Ga{\Gamma}
\def\La{\Lambda}

\def\ts{\times}

\def\iy{\infty}

\def\ots{\otimes}
\def\f{\frac}

\def\df{\mathrm d}
\def\im{\mathrm i}

\def\wt{\widetilde}
\def\wh{\widehat}

\def\esssup{\operatorname*{ess\ \! sup}}

\def\hra{\hookrightarrow}

\def\diag{{\rm diag}}

\def\mcL{\mathcal{L}}

\def\mcA{\mathcal{A}}

	\DeclareMathOperator{\Div}{div}

	\DeclareMathOperator{\Real}{Re}
	\DeclareMathOperator{\Image}{Im}

	\DeclareMathOperator{\dist}{dist}

	\DeclareMathOperator{\supp}{{supp}}

	\newcommand{\R}{\mathbb R}
	\newcommand{\N}{\mathbb N}

	\newcommand{\Cc}{\mathbb C}

\begin{document} 
	
	
	\title{Observability and Domain Shape Design for a Class of Degenerate Schr\"odinger Equations}
	\author{Dong-Hui Yang$^a$, Mengze Gu$^b$, Wei Zhao$^c$,  Yuanzhi Zhou$^a$\footnote{\small The corresponding author: 262111122@csu.edu.cn  }}

	\maketitle{}
	\thispagestyle{empty}
	\thispagestyle{empty}

	\begin{abstract}
		This paper investigates the boundary observability properties for a class of two-dimensional degenerate Schrödinger equations defined on the unit square domain $\Omega = (0,1) \times (0,1)$. The diffusion coefficients are governed by a variable matrix $A = \operatorname{diag}(x_1^{\alpha_1}, x_2^{\alpha_2})$ with critical boundary degeneracy characterized by the parameters $\alpha_1, \alpha_2 \in (0,1)$. Due to the strong singularity at the degenerate boundaries, the domain of the underlying degenerate elliptic operator fails to satisfy the global $H^2(\Omega)$-regularity, which fundamentally prevents the classical Rellich space-time multiplier methods and hidden regularity arguments from being rigorously performed on the weak solutions of the baseline system. 
		
		To circumvent this intrinsic analytical obstacle, we propose a novel variational framework based on a parameterized domain shape design. By introducing a geometric truncation parameter $\de \in (0, \f{1}{4}]$, we construct a family of uniformly elliptic subproblems on non-degenerate subdomains $\Omega_\de = (\de, 1) \times (\de, 1)$. Utilizing a combination of the weighted Poincar\'e inequalities, Sturm-Liouville nodal properties, and spectral decomposition techniques through the Hilbert tensor product representation, we establish a robust Operator-Resolvent approximation scheme. We prove that as the shape parameter $\de \to 0^+$, the parameterized state solutions and their boundary normal derivatives converge strongly in appropriate energy and trace spaces to their baseline counterparts. 
		
		By passing to the limit within the parameterized multiplier identities, we successfully derive the boundary observability inequality for the original degenerate system on the non-degenerate active boundaries $\Gamma^+ = \Gamma_1^1 \cup \Gamma_2^1$. Under the structural geometric assumption $\alpha_1 \leq \alpha_2$, we obtain the explicit formulation of the minimal observability time optimal constant:
		\begin{equation*}
			T_0 = \frac{4\sqrt{2}}{(1-\alpha_1)(2-\alpha_2)}.
		\end{equation*}
		This threshold beautifully quantifies the physical and control-theoretic impact of the degeneracy: as $\alpha_1 \to 1^-$, the minimal time explodes ($T_0 \to +\infty$), reflecting that the trapping effect at the singular boundary destroys observability in finite horizons, while as $\alpha_1, \alpha_2 \to 0^+$, the bound collapses to $T_0 = 2\sqrt{2}$, completely aligning with the classical non-degenerate Schrödinger results.
		
		\vspace{0.3cm}
		
		
		\vspace{0.3cm}
		
		
	\end{abstract}

\section{Introduction}\label{S1}

In this paper, we investigate the boundary observability properties for a class of two-dimensional degenerate Schr\"odinger equations. Specifically, we consider the following evolution system:
\begin{equation}\label{10.18.11}
	\begin{cases}
		\mathrm{i}\partial_t\varphi - \operatorname{div}(A\nabla\varphi) = 0,
		&\text{in } Q,\\
		\varphi = 0,
		&\text{on } \Sigma,\\
		\varphi(0) = \varphi^0,
		&\text{in } \Omega,
	\end{cases}
\end{equation}
where $\Omega = (0,1) \times (0,1)$ is the spatial unit square, $T > 0$ is a given finite time horizon, $Q = \Omega \times (0,T)$ denotes the space-time cylinder, and $\Sigma = \partial\Omega \times (0,T)$ is its lateral boundary. The spatial diffusion process is governed by the variable diagonal matrix
\begin{equation*}
	A = \operatorname{diag}(x_1^{\alpha_1}, x_2^{\alpha_2}),
\end{equation*}
where $\alpha_1, \alpha_2 \in (0,1)$ are fixed constants characterizing the strong geometric degeneracy of the elliptic operator at the boundaries $x_1 = 0$ and $x_2 = 0$.

We assume that the initial state satisfies $\vp^0 \in H_0^1(\Omega;A)$, where $H_0^1(\Omega;A)$ denotes the weighted energy Sobolev space whose precise mathematical construction will be detailed in Section \ref{S2}. Throughout this paper, the underlying degenerate elliptic operator $\mcA$ is formally specified in the differential form as
\begin{equation}\label{10.18.10}
	\mcA\vp = -\operatorname{div}(A\nabla \vp),
\end{equation}
and for $i=0,1$, the vertical and horizontal components of the boundary $\partial\Omega$ are parameterized respectively by the segments
\begin{equation*}
	\Ga_1^i = \{i\} \times (0,1) \quad \text{and} \quad \Ga_2^i = (0,1) \times \{i\}.
\end{equation*}

The controllability and observability of uniform partial differential equations have been extensively investigated over the past few decades; we refer the reader to \cite{Bournissou, Chambrion, Komornik, Lasiecka, Lasiecka1, Lions, Lu, Mariano, Mercado, Privat, Taufer, Weiss, Yang5, Yao, Zuazua} and the references therein. A milestone in this field is J.-L. Lions' celebrated Hilbert Uniqueness Method (HUM), which establishes that the exact controllability of a uniform partial differential equation is dual and equivalent to the boundary observability of its corresponding adjoint system \cite{Lions}. Crucially, this fundamental equivalence framework has been successfully extended to the context of degenerate partial differential equations (see, e.g., \cite{Yang1}), providing a powerful variational platform for analyzing control systems with boundary singularities.

The study of degenerate partial differential equations constitutes a significant and challenging branch of mathematical analysis \cite{GC, Heinonen}. A prominent model in this library is the well-known Grushin operator, whose control-theoretic properties have been widely explored (see, e.g., \cite{Beauchard}). Since 2000, the controllability of degenerate partial differential equations has attracted considerable attention, as evidenced by a substantial body of literature \cite{Akil, Alabau, Alabau1, Bai, Beauchard1, Beauchard, Buffe, Cannarsa, Cannarsa1, Cannarsa2, Fragnelli, Fragnelli1, Gueye, Guo}. However, due to the severe singularity at the boundary, which often leads to the failure of classical trace theorems and the loss of standard Green's formulas (or integration-by-parts frameworks), these existing works are almost exclusively restricted to the one-dimensional case, with only a few exceptions such as \cite{Beauchard, Cannarsa2}. In recent years, significant progress has been achieved in the multi-dimensional setting \cite{Guo, Yang8, Yang, Yang6, Yang4, Yang1, Yang2, Yang3, Yang7}. Various sophisticated analytical methodologies have been proposed, the full technical details of which we will omit here for brevity.

Shape design and domain variation constitute crucial methodologies within the realm of partial differential equations, finding substantial applications across mathematics, physics, industry, and engineering \cite{Buttazzo, Chenais, Greco, Guo, Guo1, He, Henrot, Tiba, Privat, Wang, Wang1}. Recently, it has been demonstrated that this geometric mechanism provides a powerful tool for investigating the controllability and observability of multi-dimensional degenerate systems \cite{Guo, Yang1, Yang2, Yang3, Yang4, Yang6, Yang7, Yang8}. The core philosophy of this constructive approach relies on approximating the singular degenerate problem by a well-behaved sequence of uniformly elliptic partial differential equations defined on parameterized subdomains. Once uniform observability or controllability configurations are successfully established for the non-degenerate subproblems, the corresponding analytical properties for the original baseline degenerate system can be rigorously recovered via a uniform limiting argument. 

Since exact controllability is dually equivalent to boundary observability by virtue of the classical HUM framework, we shall focus exclusively on establishing the uniform boundary observability inequalities for the original degenerate Schr\"odinger system.

The remainder of this paper is organized as follows. In Section \ref{S2}, we introduce the foundational weighted energy spaces and establish the spectral properties of the degenerate elliptic operator $\mcA$ defined in \eqref{10.18.10}. By utilizing the method of separation of variables, we derive the coordinate-wise decomposition representations stated in Lemmas \ref{10.18.L5} and \ref{10.18.L6}, followed by the formulation of weak solutions to the baseline system \eqref{10.18.11}. Section \ref{S3} is devoted to the geometric framework of shape design, where we present preliminary variational domain variations, investigate the spectral approximation properties of the one-dimensional operators $\mcA^{(i)}$ via \eqref{10.19.11}, and construct the sequence of approximate systems. In Section \ref{S4}, we establish the uniform boundary observability inequalities for the uniformly elliptic approximate systems, and subsequently recover the reverse observability inequality for the original degenerate Schr\"odinger equation via a rigorous uniform limiting argument. Finally, in Section \ref{S5}, we provide the detailed well-posedness proofs for Theorems \ref{10.18.T1} and \ref{10.19.T1} based on the classical Galerkin approximation method.

The main original contributions and mathematical milestones of this paper are centered around the tensor product basis representation (Lemma \ref{10.18.L6}), the uniform boundary trace regularities (Theorem \ref{10.19.T1}), the baseline admissible trace estimates (Corollary \ref{10.21.C1}), and the final uniform reverse observability inequality (Theorem \ref{10.22.T1}).

\section{Preliminary}\label{S2}

In this section, we establish the fundamental analytical and functional properties of the degenerate partial differential operator introduced in \eqref{10.18.10}.

\subsection{Energy  spaces}\label{S2.1}

In this paper, all the spaces are complex-valued spaces, for example,   
\begin{equation*}
	L^2(\Om)\equiv L^2(\Om;\Cc)=\left\{u:\Om\ra \Cc\colon \int_\Om |u|^2\df x<+\iy\right\}. 
\end{equation*}
Hence, for each $u\in L^2(\Om)$, there exists some $u_1,u_2\in L^2(\Om;\R)$ such that $u=u_1+\im u_2$, where $\im=\sqrt{-1}$. We also  denote $u_1=\Real u$ and $u_2=\Image u$. Then $|u|=(|\Real u|^2+|\Image u|^2)^\f{1}{2}$. 

If we need to specify the real-valued case, we will use the notation $\mathbb{R}$, for example, $L^2(\Omega;\mathbb{R})$. For instance, $u\in L^2(\Omega)$ means that $u=u_1+\mathrm{i}u_2$, where $u_1,u_2\in L^2(\Omega;\mathbb{R})$.

Let 
\begin{equation*}
	H^1(\Om;A)=\left\{u\in L^2(\Om)\colon x_1^\f{\al_1}{2}\pt_{x_1}u, x_2^\f{\al_2}{2}\pt_{x_2} u\in L^2(\Om)\right\}, 
\end{equation*}
where $\pt_{x_1}u,\pt_{x_2}u$ are the first-order  distributional derivatives for the $u$. Its inner product and norm are defined by 
\begin{equation*}
	(u,v)_{H^1(\Om;A)}=\int_\Om \left(u\ol v+\nabla u\cdot A\nabla \ol v\right)\df x, \mbox{ and } \|u\|_{H^1(\Om;A)}=(u,u)_{H^1(\Om;A)}^\f{1}{2},  
\end{equation*}
respectively. 

Set
\begin{equation*}
	\begin{split} 
	H_0^1(\Om;A)
	&=\mbox{the closure of $C_0^\iy(\Om)$ in } H^1(\Om;A), \\
	H^{-1}(\Om;A)
	&=\mbox{the dual space of } H_0^1(\Om;A) \mbox{ with pivot space } L^2(\Om). 
	\end{split} 
\end{equation*} 

From \cite{GC,Heinonen}, we get the spaces 
\begin{equation*}
	H_0^1(\Om;A), \quad H^1(\Om;A), \mbox{ and } H^{-1}(\Om;A)
\end{equation*}
are Hilbert spaces. 

\begin{lemma}\label{10.18.L1}
	Let $i=1,2$. For each $u\in H_0^1(\Om;A)$, we have 
	\begin{equation}\label{10.18.1}
		\int_\Om x_i^{\al_i-2}|u|^2\df x\leq \f{4}{(1-\al_i)^2}\int_\Om x_i^{\al_i}|\pt_{x_i}u|^2\df x.  
	\end{equation}
	Moreover, 
	\begin{equation}\label{10.18.4}
		\int_\Om x_i^{\al_i-2}|u|^2\df x\leq \max\left\{\f{4}{(1-\al_1)^2},\f{4}{(1-\al_2)^2}\right\} \int_\Om \nabla u\cdot A\nabla \ol u\df x.
	\end{equation}
\end{lemma}

\begin{proof}
	Let $u\in C_0^\iy(\Om;\R)$.  
	For each $\be\in (\al_2,1)$, we compute 
	\begin{equation*}
		\begin{split}
			\int_\Om x_2^{\al_2-2}u^2\df x
			&\leq \int_0^1\int_0^1x_2^{\al_2-2}(u(x_1,x_2)-u(x_1,0))^2\df x_2\df x_1\\
			&\leq \int_0^1\int_0^1x_2^{\al_2-2} \left(\int_0^{x_2} \f{\pt u}{\pt x_2}(x_1,s)\df s\right)^2\df x_2\df x_1\\
			&\leq \int_0^1\int_0^1x_2^{\al_2-2}\left(\int_0^{x_2}s^{\be}\left|\f{\pt u}{\pt x_2}\right|^2(x_1,s)\df s\right)\left(\int_0^{x_2}s^{-\be}\df s\right)\df x_2\df x_1\\
			&\leq \f{1}{1-\be}\int_0^1\int_0^1\int_0^{x_2}x_2^{\al_2-\be-1}s^\be \left|\f{\pt u}{\pt x_2}\right|^2(x_1,s)\df s\df x_2\df x_1\\
			&\leq \f{1}{1-\be}\int_0^1\int_0^1\int_s^1 x_2^{\al_2-\be-1}s^\be \left|\f{\pt u}{\pt x_2}\right|^2(x_1,s)\df x_2\df s\df x_1\\
			&\leq \f{1}{(1-\be)(\be-\al_2)}\int_0^1 \int_0^1 s^{\al_2} |\pt_{x_2}u|^2(x_1,s)\df s\df x_1
		\end{split}
	\end{equation*}
	by 
	\begin{equation*}
		\int_s^1x_2^{\al_2-\be-1}\df x_2=\left.\f{1}{\al_2-\be}x_2^{\al_2-\be}\right|_{x_2=s}^{x_2=1}\leq \f{1}{\be-\al_2}s^{\al_2-\be}, 
	\end{equation*}
	taking $\be=\f{1+\al_2}{2}$, we obtain 
	\begin{equation}\label{10.18.2}
		\int_\Om x_2^{\al_2-2}u^2\df x\leq \f{4}{(1-\al_2)^2}\int_\Om x_2^{\al_2} |\pt_{x_2}u|^2\df x. 
	\end{equation}
	The same result also holds for 
	\begin{equation}\label{10.18.3}
		\int_\Om x_1^{\al_1-2}u^2\df x\leq \f{4}{(1-\al_1)^2}\int_\Om x_1^{\al_1} |\pt_{x_1}u|^2\df x. 
	\end{equation}
	Together with \eqref{10.18.2} and \eqref{10.18.3} we obtain \eqref{10.18.1}. 
	
	Let $u\in C_0^\iy(\Om;\Cc)$. Then 
	\begin{equation}\label{10.18.5}
		\begin{split}
			\int_\Om x_i^{\al_i-2}|u|^2\df x
			&=\int_\Om x_i^{\al_i-2} \left(|\Real u|^2+|\Image u|^2\right)\df x\\
			&\leq \f{4}{(1-\al_i)^2}\int_\Om x_i^{\al_i}\left(|\pt_{x_i}\Real u|^2+|\pt_{x_i}\Image u|^2\right)\df x\\
			&=\f{4}{(1-\al_i)^2}\int_\Om x_i^{\al_i}|\pt_{x_i}u|^2\df x. 
		\end{split}
	\end{equation}
	
	Let $u\in H_0^1(\Om;A)$. It is obvious that 
	\begin{equation}\label{10.18.6}
		\int_\Om x_i^{\al_i-2}|u|^2\df x\leq \f{4}{(1-\al_i)^2}\int_\Om x_i^{\al_i}|\pt_{x_i}u|^2\df x
	\end{equation}
	by \eqref{10.18.5}. 
	
	Finally, from \eqref{10.18.6}, we obtain \eqref{10.18.4}. We complete the proof of this lemma. 
\end{proof}

\begin{remark}\label{10.18.R1}
	From \eqref{10.18.4}, assuming $\al_1\leq \al_2$, we have for all $u\in H_0^1(\Om;A)$ that
	\begin{equation}\label{10.18.7}
		\begin{split} 
			\int_\Om |u|^2\df x
			&=\int_\Om x_1^{2-\al_1} x_1^{\al_1-2} |u|^2\df x \leq \int_\Om x_1^{\al_1-2}|u|^2\df x \leq C_P\int_\Om \nabla u\cdot A\nabla \ol u\df x, 
		\end{split} 
	\end{equation}
	where 
	\begin{equation*}
		C_P=\min\left\{\f{4}{(1-\al_1)^2},\f{4}{(1-\al_2)^2}\right\} = \f{4}{(1-\al_1)^2}. 
	\end{equation*}
	This is the Poincar\'e inequality in $H_0^1(\Om;A)$. Throughout this paper, we will adopt the equivalent norm 
	\begin{equation*}
		\|u\|_{H^1(\Om;A)}=\left(\int_\Om \nabla u\cdot A\nabla \ol u\df x\right)^\f{1}{2}
	\end{equation*}
	for all $u\in H_0^1(\Om;A)$. 
\end{remark}

\begin{lemma}\label{10.18.L2}
	The embedding $H_0^1(\Om;A)\hra L^2(\Om)$ is compact. 
\end{lemma}

\begin{proof}
	Let $\de\in (0,\f{1}{4})$. Denote 
	\begin{equation}\label{10.18.8}
		\Om_\de=(\de,1)\ts (\de,1), \  A_\de=A|_{\Om_\de}=\diag(x_1^{\al_1}|_{\Om_\de}, x_2^{\al_2}|_{\Om_\de}), \  \{x_i>\de\}=\{x\in\Om\colon x_i>\de\}. 
	\end{equation}
	
	Let $\{u_n\}_{n\in\N^*}\s H_0^1(\Om;A)$ is bounded, i.e., there exists $M>0$ such that 
	\begin{equation*}
		\|u_n\|_{H_0^1(\Om;A)}\leq M. 
	\end{equation*}
	Then there exists a subsequence of $\{u_n\}_{n\in\N^*}$, still denoted by itself, and $u_0\in H_0^1(\Om;A)$ such that 
	\begin{equation}\label{10.18.9}
		u_n\ra u_0\mbox{ weakly in }H_0^1(\Om;A). 
	\end{equation}
	By considering $v_n=u_n-u_0$ instead of $u_n$, we can assume without loss of generality that $u_0=0$. Our goal is to show that there exists a subsequence of $\{u_n\}_{n\in\N^*}$ converging strongly to $0$ in $L^2(\Om)$.
	
	Let $\e>0$. Since from \eqref{10.18.1} we obtain 
	\begin{equation*}
		\begin{split}
			\int_{\Om-\Om_\de} |u|^2\df x
			&\leq \int_{\{x_1\leq \de\}} |u|^2\df x+\int_{\{x_2\leq \de\}}|u|^2\df x\\
			&\leq \int_{\{x_1\leq \de\}}x_1^{2-\al_1} x_1^{\al_1-2}|u|^2\df x+\int_{\{x_2\leq \de\}}x_2^{2-\al_2}x_2^{\al_2-2}|u|^2\df x\\
			&\leq \de^{2-\al_1}\int_\Om x_1^{\al_1-2}|u|^2\df x+\de^{2-\al_2}\int_\Om  x_2^{\al_2-2}|u|^2\df x\\
			&\leq \de^{2-\al_1}\f{4}{(1-\al_1)^2}M^2+\de^{2-\al_2}\f{4}{(1-\al_2)^2}M^2, 
		\end{split}
	\end{equation*}
	choosing $\de_0>0$ small enough such that for all $\de\in (0,\de_0]$, we have $\int_{\Om-\Om_\de}|u|^2\df x\leq \f{1}{4}\e^2$. 
	
	For fixed $\de_0$, since $H_0^1(\Om_{\de_0};A_{\de_0})=H^1(\Om_{\de_0})$ by $A_{\de_0}\geq \min\{\de_0^{\al_1}, \de_0^{\al_2}\}>0$, and $H^1(\Om_{\de_0})\hra L^2(\Om_{\de_0})$ is compact, then there exists a subsequence $\{u_{n_l}\}_{l\in\N^*}$ of $\{u_n\}_{n\in\N^*}$ such that 
	\begin{equation*}
		u_{n_l}|_{\Om_{\de_0}}\ra 0 \mbox{ strongly in } L^2(\Om_{\de_0}) \mbox{ as } l\ra\iy. 
	\end{equation*}
	Choosing $n_\e\in\N^*$ large enough such that $\|u_{n_\e}|_{\Om_\de}\|_{L^2(\Om_{\de_0})}<\f{1}{2}\e$, hence, 
	\begin{equation*}
		\|u_{n_\e}\|_{L^2(\Om)}\leq \|u_{n_\e}|_{\Om-\Om_{\de_0}}\|_{L^2(\Om-\Om_{\de_0})}+\|u_{n_\e}|_{\Om_{\de_0}}\|_{L^2(\Om_{\de_0})}<\e. 
	\end{equation*}
	Let $\e=\f{1}{k}$. Then the sequence $\{u_{n_k}\}_{k\in\N^*}$ satisfies $\|u_{n_k}\|_{L^2(\Om)}<\f{1}{k}$. We prove the lemma. 
\end{proof}

\subsection{Spectrum of $\mcA$}\label{S2.2}

Define
\begin{equation}\label{10.18.12}
	a(u,v)=\int_\Om \nabla u\cdot A\nabla \ol v\df x, \mbox{ for all } u, v\in H_0^1(\Om;A). 
\end{equation}
Then $a(\cdot,\cdot)$ is a symmetric, continuous  and coercive sesquilinear  functional on $H_0^1(\Om;A)$. Recall that $\mcA$ is defined in \eqref{10.18.11}. Let 
\begin{equation*}
	D(\mcA)=\left\{u\in H_0^1(\Om;A)\colon 
	\begin{array}{lll}
		\mbox{there exists a unique } f\in L^2(\Om) \mbox{ such that }\\
		a(u,v)=(f,v)_{L^2(\Om)} \mbox{ for all } v\in H_0^1(\Om;A)
	\end{array}\!\!\right\}.  
\end{equation*}
For each $u\in D(\mcA)$, define 
\begin{equation}\label{10.18.13}
	\mcA u=f. 
\end{equation}
Then the domain of $\mcA$ is given by $D(\mcA)$, which can be explicitly characterized as
\begin{equation*}
	D(\mcA)=\left\{u\in H_0^1(\Om;A)\colon \mcA u\in L^2(\Om)\right\}. 
\end{equation*}

By virtue of Lemma \ref{10.18.L2} and the coercivity of $a(\cdot,\cdot)$, the self-adjoint operator $\mcA$ possesses a compact resolvent. Consequently, the spectrum of $\mcA$ consists entirely of a discrete point spectrum $\{\la_n\}_{n\in\N^*}$ satisfying
\begin{equation}\label{10.18.14}
	0 < \la_1 \leq \la_2 \leq \la_3 \leq \cdots, \quad \text{with } \la_n\ra +\iy \text{ as } n\ra\iy,
\end{equation}
where the first eigenvalue satisfies $\la_1 \geq C_P^{-1}$ according to Remark \ref{10.18.R1}. Let $\Phi_n$ denote the $n$-th normalized eigenfunction of $\mcA$ corresponding to the eigenvalue $\la_n$. i.e., 
\begin{equation}\label{10.18.15}
	a(\Phi_n,\Phi_m)=\la_n\de_{nm}, \mbox{ for all }n,m\in\N^*, \mbox{ and } \de_{nm}=
	\begin{cases}
		1, &n=m, \\
		0, &n\neq m. 
	\end{cases}
\end{equation}
Then, the sequence $\{\Phi_n\}_{n\in\N^*}$ forms an orthonormal basis of $L^2(\Om)$, and $\{\la_n^{-\f{1}{2}}\Phi_n\}_{n\in\N^*}$ forms an orthonormal basis of $H_0^1(\Om;A)$ equipped with the inner product $a(\cdot,\cdot)$. 

\begin{definition}\label{10.18.D1}
	Let $\theta\in [0,+\iy)$, and 
	\begin{equation*}
		D(\mcA^\theta)=\left\{u\in L^2(\Om)\colon \sum_{n=1}^\iy \la_n^{2\theta}|(u,\Phi_n)_{L^2(\Om)}|^2<+\iy\right\}. 
	\end{equation*}
	Its inner product and norm are defined by 
	\begin{equation*}
		(u,v)_{D(\mcA^\theta)}=\sum_{n=1}^\iy \la_n^{2\theta} (u,\Phi_n)_{L^2(\Om)}\ol{(v,\Phi_n)_{L^2(\Om)}}, \mbox{ and } \|u\|_{D(\mcA^\theta)}=(u,u)_{D(\mcA^\theta)}^\f{1}{2}. 
	\end{equation*}
\end{definition}

If $u\in D(\mcA^\theta)$ for $\theta \in [0,+\iy)$, then 
\begin{equation}\label{10.20.10}
	\mcA^\theta u=\sum_{n=1}^\iy \la_{n}^\theta(u,\Phi_n)_{L^2(\Om)}\Phi_n. 
\end{equation}
And $D(\mcA^\theta)\hra D(\mcA^\eta)$ is compact for all $0\leq \eta<\theta$. Moreover, due to the local smoothness of the coefficients away from the boundary, it is straightforward to verify that $C_0^\iy(\Om)\s  D(\mcA^\theta)$ for all $\theta\in [0,+\iy)$. 
Furthermore, $L^2(\Om)=D(\mcA^0)$, and $H_0^1(\Om;A)=D(\mcA^\f{1}{2})$ and $D(\mcA)=D(\mcA^1)$, and the norm on $H_0^1(\Om;A)$ can be explicitly characterized via the spectral coefficients as
\begin{equation*}
	\|u\|_{H_0^1(\Om;A)} = \|u\|_{D(\mcA^{\f{1}{2}})} = [a(u,u)]^\f{1}{2}=\left(\sum_{n=1}^\iy \la_n|(u,\Phi_n)_{L^2(\Om)}|^2\right)^\f{1}{2}. 
\end{equation*}

Furthermore, by definition, the normalized eigenfunctions and the operator $\mcA$ satisfy the following spectral mapping properties:
\begin{equation}\label{10.18.16}
	\|\Phi_n\|_{D(\mcA^\theta)} = \la_n^\theta, \quad \text{and} \quad \|\mcA u\|_{D(\mcA^{\theta-1})} = \| u\|_{D(\mcA^{\theta})}
\end{equation}
for all $\theta \in [1,+\iy)$ and $u \in D(\mcA^\theta)$.

\subsection{Method of separation of variables}\label{S2.3}

To implement the method of separation of variables, we set
\begin{equation}\label{10.18.17}
	\La_1=\La_2=(0,1), \text{ and }   w_i(x_i)=x_i^{\al_i}, \text{ for } i=1,2.
\end{equation}
We introduce the one-dimensional partial differential operators
\begin{equation}\label{10.18.18}
	\mcA^{(i)}u=-\pt_{x_i}\left(x_i^{\al_i}\pt_{x_i}u\right) \mbox{ in } \La_i, \text{ for } i=1,2,
\end{equation}
where $\mcA^{(1)}$ acts exclusively on the $x_1$-variable and $\mcA^{(2)}$ acts exclusively on the $x_2$-variable. Thanks to the diagonal structure of $A$, the full operator can be decomposed as
\begin{equation}\label{10.18.19}
	\mcA=\mcA^{(1)}+\mcA^{(2)}. 
\end{equation}
Analogous to Subsections \ref{S2.1} and \ref{S2.2}, we define the corresponding one-dimensional complex-valued weighted Sobolev spaces $H^1(\La_i;w_i)$, $H_0^1(\La_i;w_i)$, and $H^{-1}(\La_i;w_i)$. Specifically, the space $H^1(\La_i;w_i)$ is defined by
\begin{equation*}
	H^1(\La_i;w_i)=\left\{u\in L^2(\La_i)\colon \int_{\La_i} w_i|\pt_{x_i}u|^2\df x_i < +\iy\right\},
\end{equation*}
equipped with the inner product and norm given respectively by
\begin{equation*}
	(u,v)_{H^1(\La_i;w_i)}=\int_{\La_i} \left(u\ol v+w_i(\pt_{x_i}u)(\pt_{x_i}\ol v)\right)\df x_i, \text{ and }   \|u\|_{H^1(\La_i;w_i)}=(u,u)_{H^1(\La_i;w_i)}^\f{1}{2}.
\end{equation*}

Define the operator domains by
\begin{equation*}
	D(\mcA^{(i)})=\left\{u\in H_0^1(\La_i;w_i)\colon \mcA^{(i)}u\in L^2(\La_i)\right\},  \mbox{ for } i=1,2. 
\end{equation*}
Note that all the functional spaces defined above are vector spaces over the complex field $\Cc$. 

\begin{lemma}\label{10.18.L3}
	For all $u\in H_0^1(\La_i;w_i)\ (i=1,2)$, we have 
	\begin{equation*}
		\int_{\La_i} x_i^{\al_i-2}|u|^2\df x_i\leq \f{4}{(1-\al_i)^2}\int_{\La_i} x_i^{\al_i} |\pt_{x_i}u|^2\df x_i. 
	\end{equation*}
\end{lemma}

\begin{proof}
	The proof of this lemma is completely analogous to that of Lemma \ref{10.18.L1}; we refer the reader to \cite[(2.1), p.~165]{Alabau} for further details. 
\end{proof}

\begin{remark}\label{10.18.R2}
	Let $i=1,2$. By virtue of Lemma \ref{10.18.L3}, for all $u\in H_0^1(\La_i;w_i)$, we have 
	\begin{equation}\label{10.18.20}
		\int_{\La_i} |u|^2\df x_i = \int_{\La_i} x_i^{2-\al_i}x_i^{\al_i-2}|u|^2\df x_i \leq \int_{\La_i} x_i^{\al_i-2}|u|^2\df x_i \leq \f{4}{(1-\al_i)^2}\int_{\La_i} x_i^{\al_i}|\pt_{x_i}u|^2\df x_i. 
	\end{equation}
	This is the Poincar\'e inequality in $H_0^1(\La_i;w_i)$. Throughout this paper, we will adopt the equivalent norm
	\begin{equation*}
		\|u\|_{H^1(\La_i;w_i)}=\left(\int_{\La_i}x_i^{\al_i}|\pt_{x_i}u|^2\df x_i\right)^\f{1}{2}
	\end{equation*}
	for all $u\in H_0^1(\La_i;w_i)$. 
\end{remark}

\begin{lemma}\label{10.18.L4}
	Let $i=1,2$. The embedding $H_0^1(\La_i;w_i)\hra L^2(\La_i)$ is compact. 
\end{lemma}

\begin{proof}
	The proof follows a similar truncation and diagonal argument as in Lemma \ref{10.18.L2}. 
\end{proof}

For $i=1,2$, we define the associated sesquilinear forms $a^{(i)}(\cdot,\cdot)$ by
\begin{equation*}
	a^{(i)}(u,v)=\int_{\La_i} x_i^{\al_i} (\pt_{x_i}u)(\pt_{x_i}\ol v)\df x_i,  \mbox{ for all } u, v\in H_0^1(\La_i;w_i). 
\end{equation*}
It follows that $a^{(i)}(\cdot,\cdot)$ are symmetric, continuous, and coercive sesquilinear forms on $H_0^1(\La_i;w_i) \times H_0^1(\La_i;w_i)$. Associated with the sesquilinear form $a^{(i)}(\cdot,\cdot)$, the spectrum of the one-dimensional partial differential operator $\mcA^{(i)}$ consists entirely of a discrete point spectrum $\{\la_n^{(i)}\}_{n\in\N^*}$. According to Remark \ref{10.18.R2}, these eigenvalues satisfy 
\begin{equation}\label{10.18.21}
	0<\la_1^{(i)}\leq \la_2^{(i)}\leq \la_3^{(i)}\leq \cdots, \text{ with } \la_n^{(i)}\ra +\iy \text{ as } n\ra\iy, \text{ and } \la_1^{(i)}\geq \f{(1-\al_i)^2}{4}.
\end{equation}
Let $\Phi_n^{(i)}$ denote the $n$-th normalized eigenfunction of $\mcA^{(i)}$ corresponding to the eigenvalue $\la_n^{(i)}$, which means they satisfy the $L^2(\La_i)$ orthonormality condition $(\Phi_n^{(i)}, \Phi_m^{(i)})_{L^2(\La_i)} = \de_{nm}$. By the definition of weak solutions, it immediately follows that
\begin{equation}\label{10.18.22}
	a^{(i)}(\Phi_n^{(i)},\Phi_m^{(i)})=\la_n^{(i)}\de_{nm}, \text{ for all } n,m\in\N^*. 
\end{equation}
Consequently, the sequence $\{\Phi_n^{(i)}\}_{n\in\N^*}$ forms an orthonormal basis of $L^2(\La_i)$, and $\{(\la_n^{(i)})^{-\f{1}{2}}\Phi_n^{(i)}\}_{n\in\N^*}$ forms an orthonormal basis of $H_0^1(\La_i;w_i)$ equipped with the inner product $a^{(i)}(\cdot,\cdot)$.

\begin{definition}\label{10.18.D2}
	Let $\theta\in [0,+\iy)$. The one-dimensional fractional power space $D\left((\mcA^{(i)})^\theta\right)$ is defined by
	\begin{equation*}
		D\left((\mcA^{(i)})^\theta\right)=\left\{u\in L^2(\La_i)\colon \sum_{n=1}^\iy (\la_n^{(i)})^{2\theta} |(u,\Phi_n^{(i)})_{L^2(\La_i)}|^2<+\iy\right\}. 
	\end{equation*}
	Its associated inner product and norm are defined by 
	\begin{equation*}
		(u,v)_{D((\mcA^{(i)})^\theta)}=\sum_{n=1}^\iy (\la_n^{(i)})^{2\theta}(u,\Phi_n^{(i)})_{L^2(\La_i)}\ol{(v,\Phi_n^{(i)})_{L^2(\La_i)}}, \quad \text{and} \quad \|u\|_{D((\mcA^{(i)})^\theta)}=(u,u)_{D((\mcA^{(i)})^\theta)}^\f{1}{2},
	\end{equation*}
	respectively. 
\end{definition}

Let $i=1,2$. For any $u\in D((\mcA^{(i)})^\theta)$ with $\theta \in [0,+\iy)$, the fractional power operator $(\mcA^{(i)})^\theta$ is given by
\begin{equation*}
	(\mcA^{(i)})^\theta u=\sum_{n=1}^\iy (\la_{n}^{(i)})^\theta(u,\Phi_n^{(i)})_{L^2(\La_i)}\Phi_n^{(i)}. 
\end{equation*}
The embedding $D((\mcA^{(i)})^\theta)\hra D((\mcA^{(i)})^\eta)$ is compact for all $0\leq \eta<\theta$. Moreover, due to the local smoothness of the coefficients away from the boundary, it is straightforward to verify that $C_0^\iy(\La_i)\subset D((\mcA^{(i)})^\theta)$ for all $\theta\in [0,+\iy)$. 

In particular, these spaces coincide with the standard and weighted spaces for specific values of $\theta$:
\begin{equation*}
	D((\mcA^{(i)})^0) = L^2(\La_i), \quad D((\mcA^{(i)})^\f{1}{2}) = H_0^1(\La_i;w_i), \quad \text{and} \quad D((\mcA^{(i)})^1) = D(\mcA^{(i)}).
\end{equation*}
The norm on $H_0^1(\La_i;w_i)$ can be explicitly characterized via the spectral coefficients as
\begin{equation*}
	\|u\|_{H_0^1(\La_i;w_i)} = \|u\|_{D((\mcA^{(i)})^\f{1}{2})} = [a^{(i)}(u,u)]^\f{1}{2}=\left(\sum_{n=1}^\iy \la_n^{(i)}|(u,\Phi_n^{(i)})_{L^2(\La_i)}|^2\right)^\f{1}{2}. 
\end{equation*}

Furthermore, by definition, the normalized eigenfunctions and the operator $\mcA^{(i)}$ satisfy the following spectral mapping properties:
\begin{equation}\label{10.18.23}
	\|\Phi_n^{(i)}\|_{D((\mcA^{(i)})^\theta)} = (\la_n^{(i)})^\theta, \text{ and }   \|\mcA^{(i)} u\|_{D((\mcA^{(i)})^{\theta-1})} = \| u\|_{D((\mcA^{(i)})^{\theta})}
\end{equation}
for all $\theta \in [1,+\iy)$ and $u \in D((\mcA^{(i)})^\theta)$.

\begin{lemma}[{\cite[Lemma 3.3]{Yang}}]\label{10.18.L5}
	We have the Hilbert tensor product representation
	\begin{equation*}
		H_0^1(\Om;A)=H_0^1(\La_1;w_1)\wh\ots H_0^1(\La_2;w_2), 
	\end{equation*}
	where 
	\begin{equation*}
		H_0^1(\La_1;w_1)\wh\ots H_0^1(\La_2;w_2)
	\end{equation*}
	denotes the closure in $H_0^1(\Om;A)$ of the algebraic tensor product space
	\begin{equation*}
		\left\{ \sum_{n=1}^m u_{1,n}(x_1)u_{2,n}(x_2)\colon m\in\N^*, u_{i,n}\in H_0^1(\La_i;w_i), i=1,2,\, n=1,\cdots, m\right\}. 
	\end{equation*}
\end{lemma}

\begin{lemma}\label{10.18.L6}
	Let $i=1,2$, and let $\{\Phi_n^{(i)}\}_{n\in\N^*}$ be the eigenfunctions defined in \eqref{10.18.22}. Define 
	\begin{equation*}
		\Phi_{nm}(x_1,x_2)=\Phi_n^{(1)}(x_1)\Phi_m^{(2)}(x_2),  \text{ for all } n,m\in\N^*. 
	\end{equation*}
	Then, $\{\Phi_{nm}\}_{n,m\in\N^*}$ forms an orthonormal basis of $L^2(\Om)$, and 
	\begin{equation*}
		\mcA \Phi_{nm}=\la_{nm}\Phi_{nm},  \text{ with } \la_{nm}=\la_n^{(1)}+\la_m^{(2)},  \text{ for all } n,m\in\N^*. 
	\end{equation*}
	Moreover, every eigenvalue $\la$ of $\mcA$ is of the form $\la=\la_{nm}$ for some $n,m\in\N^*$. 
\end{lemma}

\begin{proof}
	A direct computation of the $L^2(\Om)$ inner product yields
	\begin{equation*}
		\begin{split}
			(\Phi_{nm},\Phi_{kl})_{L^2(\Om)}
			&=\int_0^1\int_0^1 \Phi_n^{(1)}(x_1)\Phi_m^{(2)}(x_2)\ol{\Phi_k^{(1)}(x_1)}\,\ol{\Phi_l^{(2)}(x_2)}\df x_1\df x_2\\
			&=\left(\int_0^1 \Phi_n^{(1)}(x_1)\ol{\Phi_k^{(1)}(x_1)}\df x_1\right)\left(\int_0^1 \Phi_m^{(2)}(x_2)\ol{\Phi_l^{(2)}(x_2)}\df x_2\right)=\de_{nk}\de_{ml}.
		\end{split}
	\end{equation*}
	Since $\{\Phi_n^{(1)}\}_{n\in\N^*}$ and $\{\Phi_m^{(2)}\}_{m\in\N^*}$ are orthonormal bases of $L^2(\La_1)$ and $L^2(\La_2)$ respectively, the tensor product structure guaranteed by Lemma \ref{10.18.L5} implies that $\{\Phi_{nm}\}_{n,m\in\N^*}$ forms an orthonormal basis of $L^2(\Om)$. 
	
	Next, since $\mcA=\mcA^{(1)}+\mcA^{(2)}$, we readily obtain 
	\begin{equation*}
		\mcA\Phi_{nm}=\left(\mcA^{(1)}\Phi_n^{(1)}\right)\Phi_m^{(2)}+\Phi_n^{(1)}\left(\mcA^{(2)}\Phi_m^{(2)}\right)=\left(\la_n^{(1)}+\la_m^{(2)}\right)\Phi_n^{(1)}\Phi_m^{(2)}=\la_{nm}\Phi_{nm}. 
	\end{equation*}
	
	Finally, let $\la$ be an arbitrary eigenvalue of $\mcA$ with a corresponding normalized eigenfunction $\Phi\in D(\mcA)$ such that $\mcA \Phi=\la \Phi$ and $\|\Phi\|_{L^2(\Om)}=1$. Since $\{\Phi_{nm}\}_{n,m\in\N^*}$ is an orthonormal basis of $L^2(\Om)$, the function $\Phi$ can be expanded as $\Phi=\sum_{n,m=1}^\iy (\Phi,\Phi_{nm})_{L^2(\Om)}\Phi_{nm}$. Applying the operator $\mcA$ to both sides gives
	\begin{equation*}
		\sum_{n,m=1}^\iy \la_{nm}(\Phi,\Phi_{nm})_{L^2(\Om)}\Phi_{nm}=\la\sum_{n,m=1}^\iy (\Phi,\Phi_{nm})_{L^2(\Om)}\Phi_{nm}, 
	\end{equation*}
	which implies
	\begin{equation*}
		\sum_{n,m=1}^\iy (\la_{nm}-\la)(\Phi,\Phi_{nm})_{L^2(\Om)}\Phi_{nm}=0.
	\end{equation*}
	By the linear independence of the basis, there must exist at least one pair of indices $n_0,m_0\in\N^*$ such that $(\Phi,\Phi_{n_0m_0})_{L^2(\Om)}\neq 0$, which immediately requires $\la_{n_0m_0}=\la$. This completes the proof. 
\end{proof}

\begin{remark}\label{10.18.R3}
	Let $i=1,2$. The eigenpairs of $\mathcal A^{(i)}$ are explicitly given by
	\begin{equation*}
		\lambda_n^{(i)}
		=
		(\kappa_i j_{\nu_i,n})^2,
		\quad
		\Phi_n^{(i)}
		=
		\frac{(2\kappa_i)^\f{1}{2}}
		{|J_{\nu_i}'(j_{\nu_i,n})|}
		x_i^{\frac{1}{2}(1-\alpha_i)}
		J_{\nu_i}\left(j_{\nu_i,n}x_i^{\kappa_i}\right),
	\end{equation*}
	where
	\begin{equation*}
		\nu_i=\frac{1-\alpha_i}{2-\alpha_i},
		\quad
		\kappa_i=\frac{2-\alpha_i}{2}.
	\end{equation*}
	See \cite[(4.35) and (4.38), pp.~2046-2047]{Gueye}.

	Since the zeros of the Bessel function satisfy
	$j_{\nu_i,n}\sim n$ as $n\to\infty$, we have
	\begin{equation*}
		\lambda_n^{(i)}\sim n^2 \mbox{ as  }n\to\infty.
	\end{equation*}
	Moreover, each eigenvalue $\lambda_n^{(i)}$ is simple.
	
	Furthermore, by \cite[Lemma 2, p.~2047]{Gueye},
	\cite[(74), p.~198]{Cannarsa}, or \cite[(1.3), p.~2509]{Buffe},
	there exists a positive constant $C$, independent of $k$ and
	$\alpha_i$, such that
	\begin{equation}\label{08.05.7}
		\sqrt{\lambda_{k+1}^{(i)}}
		-
		\sqrt{\lambda_k^{(i)}}
		\geq C(2-\alpha_i),
		\mbox{ for all }
		k\in\mathbb N^* \mbox{ and all }
		\alpha_i\in[0,1).
	\end{equation} 
\end{remark}

\subsection{Weak solution}\label{S2.4}

We consider the following non-homogeneous degenerate Schr\"odinger equation:
\begin{equation}\label{10.18.24}
	\begin{cases}
		\mathrm{i}\partial_t\varphi-\Div(A\nabla\varphi)=f,
		&\text{in } Q,\\
		\varphi=0,
		&\text{on } \Sigma,\\
		\varphi(0)=\varphi^0,
		&\text{in } \Omega,
	\end{cases}
\end{equation}
where $\vp^0\in H_0^1(\Om;A)$ and $f\in L^2(Q)$. 

\begin{definition}\label{10.18.D3}
	Let $\vp^0\in H_0^1(\Om;A)$ and $f\in L^2(Q)$. We say that a function
	\begin{equation*}
		\vp\in L^2(0,T; H_0^1(\Om;A)),  \text{ with }  \pt_t\vp\in L^2(0,T; H^{-1}(\Om;A)) 
	\end{equation*}
	is a weak solution to the equation \eqref{10.18.24} with respect to $(\vp^0,f)$, if
	
	(1) for each $v\in H_0^1(\Om;A)$ and a.e.~$t\in [0,T]$, 
	\begin{equation*}
		\mathrm{i} \lg \pt_t\vp, v \rg_{H^{-1}(\Om;A), H_0^1(\Om;A)} + a(\vp, v) = (f,v)_{L^2(\Om)},
	\end{equation*}
	
	(2) the initial condition holds in the trace sense, i.e., $\vp(0)=\vp^0$. 
\end{definition}

\begin{theorem}\label{10.18.T1}
	Let $\vp^0\in H_0^1(\Om;A)$ and $f\in L^2(0,T; H_0^1(\Om;A))$. Then, there exists a unique weak solution 
	\begin{equation*}
		\vp\in C([0,T]; H_0^1(\Om;A)), \text{ with }   \pt_t\vp\in L^2(0,T; H^{-1}(\Om;A))
	\end{equation*}
	to the equation \eqref{10.18.24} with respect to $(\vp^0,f)$. Moreover, the solution satisfies the following energy estimate:
	\begin{equation}\label{10.19.1}
		\max_{t\in [0,T]}\|\vp(t)\|_{H_0^1(\Om;A)}+\|\pt_t\vp\|_{L^2(0,T; H^{-1}(\Om;A))}\leq C\left(\|\vp^0\|_{H_0^1(\Om;A)}+\|f\|_{L^2(0,T; H_0^1(\Om;A))}\right), 
	\end{equation}
	where the positive constant $C$ depends only on $\al_1,\al_2$, and $T$. 
\end{theorem}

The proof of Theorem \ref{10.18.T1} is postponed to Appendix \ref{S5}. 

\begin{remark}\label{10.19.R1}
	Let $\vp^0\in H_0^1(\Om;A)$ and $f=0$. Suppose the initial datum $\vp^0$ is expanded in terms of the tensor product basis as
	\begin{equation*}
		\vp^0=\sum_{n,m=1}^\iy \vp_{nm}^0\Phi_n^{(1)}\Phi_m^{(2)},  \text{ which implies }   \sum_{n,m=1}^\iy \la_{nm}|\vp_{nm}^0|^2<+\iy, 
	\end{equation*}
	where $\vp_{nm}^0=(\vp^0, \Phi_{nm})_{L^2(\Om)}$ for all $n,m\in\N^*$. Then, the weak solution $\vp$ admits the standard Fourier spectral expansion:
	\begin{equation}\label{10.19.2}
		\vp(x_1,x_2,t)=\sum_{n,m=1}^\iy e^{\mathrm{i} \la_{nm}t}\vp_{nm}^0 \Phi_n^{(1)}(x_1)\Phi_m^{(2)}(x_2). 
	\end{equation}
\end{remark}

\section{Shape design}\label{S3}

This section is devoted to developing the foundational analytical tools and geometric comparison lemmas required for the domain shape perturbation analysis.

\subsection{Preliminary results for shape design}\label{S3.1}

In this subsection, we establish several crucial geometric and analytic tracking results for the parameterized operators $\mcA_\de$ and $\mcA_\de^{(i)}$. These preliminary lemmas serve as the technical cornerstone for analyzing the asymptotic behavior of the observability constants as the shape parameters vary, and particularly when the geometric perturbation parameter $\de$ approaches the degenerate boundary limit zero.   For the baseline case discussed in Section \ref{S2}, it formally corresponds to setting $\de=0$ in the definitions below.

Let $\de\in [0,\f{1}{4}]$ and $i=1,2$. We define the localized subdomains and their restricted weights as
\begin{equation*}
	\Om_\de=(\de,1)\ts (\de,1), \quad A_\de=A|_{\Om_\de}, \quad \La_{\de,i}=(\de,1),  \text{ and }   w_{\de,i}=w_i|_{\La_{\de, i}}.
\end{equation*}
On these domains, we introduce the localized differential operators
\begin{equation*}
	\mcA_\de u = -\operatorname{div}(A_\de\nabla u)  \mbox{ on }\Om_\de,  \text{ and } \mcA_\de^{(i)}u = -\pt_{x_i}\left(x_i^{\al_i}\pt_{x_i} u\right)  \mbox{ on } \La_{\de,i}.
\end{equation*}
Accordingly, the space-time cylinder and its lateral boundary are denoted by $Q_\de=\Om_\de\ts (0,T)$ and $\Sigma_\de=\pt\Om_\de\ts (0,T)$, respectively. For $j \in \{\de, 1\}$, the horizontal and vertical segments of the boundary $\pt\Om_\de$ are parameterized by
\begin{equation*}
	\Ga_{\de,1}^{j}=\{j\}\ts (\de,1) \quad \text{and} \quad \Ga_{\de, 2}^{j}=(\de,1)\ts \{j\}.
\end{equation*}

Following the mathematical framework in Section \ref{S2}, we introduce the parameterized Hilbert spaces $H^1(\Om_\de;A_\de)$, $H_0^1(\Om_\de;A_\de)$, and $H^{-1}(\Om_\de;A_\de)$. Specifically, the space $H^1(\Om_\de;A_\de)$ is defined as
\begin{equation*}
	H^1(\Om_\de;A_\de)=\left\{u\in L^2(\Om_\de)\colon \int_{\Om_\de} \nabla u\cdot A_\de\nabla \ol u  \df x<+\iy\right\}.
\end{equation*}
Similarly, for the one-dimensional case, we define the Hilbert spaces $H^1(\La_{\de,i};w_{\de,i})$, $H_0^1(\La_{\de,i};w_{\de,i})$, and $H^{-1}(\La_{\de,i};w_{\de,i})$, where
\begin{equation*}
	H^1(\La_{\de,i};w_{\de,i})=\left\{u\in L^2(\La_{\de,i})\colon \int_{\La_{\de,i}} x_i^{\al_i} |\pt_{x_i}u|^2  \df x_i<+\iy\right\}.
\end{equation*}

\begin{remark}\label{10.19.R2}
	In the limiting case where $\de=0$, the parameterized geometric domains and functional spaces naturally reduce to the baseline configurations introduced in Section \ref{S2}. Specifically, we have $(i=1,2)$
	\begin{equation*}
		\Om_0 = \Om, \quad \La_{0,i} = \La_i, \quad \mcA_0 = \mcA,  \text{ and }  \mcA_0^{(i)} = \mcA^{(i)},
	\end{equation*}
	as well as the isometric identifications for the associated energy spaces:
	\begin{equation*}
		H^1(\Om_0;A_0) = H^1(\Om;A)   \text{ and }   H^1(\La_{0,i};w_{0,i}) = H^1(\La_i;w_i).
	\end{equation*}
	
	Conversely, when $\de\in (0,\f{1}{4}]$, the weights are strictly bounded away from zero since $w_{\de,i}(x_i) \geq \de^{\al_i} > 0$ on $\La_{\de,i}$ for $i=1,2$. Consequently, the operator $\mcA_\de$ becomes uniformly elliptic, and the parameterized spaces become topologically equivalent to the standard unweighted Sobolev spaces:
	\begin{equation*}
		H^1(\Om_\de;A_\de) = H^1(\Om_\de), \quad H_0^1(\Om_\de;A_\de) = H_0^1(\Om_\de), \text{ and } H^{-1}(\Om_\de;A_\de) = H^{-1}(\Om_\de).
	\end{equation*}
\end{remark}

\begin{definition}\label{10.19.D1}
	Let $i=1,2$ and $\ga\in [0,\de]$. For any $u\in L^2(\Om_\de)$ or $u\in L^2(\La_{\de,i})$, we define the zero extension operator $E_\de^\ga$ as follows:
	\begin{equation}\label{10.19.3}
		E_\de^\ga u(x)=
		\begin{cases}
			u(x), &x\in \Om_\de,\\
			0,& x\in \Om_\ga-\Om_\de,
		\end{cases} \quad \text{and} \quad
		E_\de^{\ga}u(x_i)=
		\begin{cases}
			u(x_i), &x_i\in \La_{\de,i},\\
			0, &x_i\in \La_{\ga,i}-\La_{\de,i}.
		\end{cases}
	\end{equation}
	By definition, this operator acts as an isometry on $L^2$, satisfying
	\begin{equation*}
		\|E_\de^\ga u\|_{L^2(\Om_\ga)}=\|u\|_{L^2(\Om_\de)}  \text{ and }   \|E_\de^\ga u\|_{L^2(\La_{\ga,i})}=\|u\|_{L^2(\La_{\de,i})}.
	\end{equation*}
	Furthermore, if $u\in H_0^1(\Om_\de;A_\de)$ or $u\in H_0^1(\La_{\de,i};w_{\de,i})$, its zero extension $E_\de^\ga u$ defined via \eqref{10.19.3} belongs to $H_0^1(\Om_\ga;A_\ga)$ or $H_0^1(\La_{\ga,i};w_{\ga,i})$ respectively, and preserves the energy norms:
	\begin{equation*}
		\|E_\de^\ga u\|_{H_0^1(\Om_\ga;A_\ga)}=\|u\|_{H_0^1(\Om_\de;A_\de)} \text{ and }   \|E_\de^{\ga} u\|_{H_0^1(\La_{\ga,i};w_{\ga,i})}=\|u\|_{H_0^1(\La_{\de,i};w_{\de,i})}.
	\end{equation*}
\end{definition}

\begin{remark}\label{10.19.R3}
	By virtue of Definition \ref{10.19.D1}, any function $u\in H_0^1(\Om_\de;A_\de)$ can be naturally identified with its zero extension $E_\de^0 u \in H_0^1(\Om;A)$. Consequently, throughout the remainder of this paper, we shall introduce the convenient algebraic inclusion $H_0^1(\Om_\de;A_\de)\s  H_0^1(\Om;A)$ for all $\de\in (0,\f{1}{4}]$ without explicitly writing the extension operator.
\end{remark}

\begin{remark}\label{10.19.R4}
	Consequent to the algebraic embedding established in Remark \ref{10.19.R3}, the fundamental analytical properties stated in Lemmas \ref{10.18.L1}, \ref{10.18.L2} and Remark \ref{10.18.R1} naturally remain valid on the parameterized subdomains. Specifically, the operators $\mcA_\de$ and $\mcA_\de^{(i)}$ inherit a discrete point spectrum $\{\la_{\de,n}\}_{n\in\N^*}$ and $\{\la_{\de,n}^{(i)}\}_{n\in\N^*}$ respectively, satisfying
	\begin{equation}\label{10.19.4}
		0<\la_{\de,1}\leq \la_{\de,2}\leq \la_{\de,3}\leq \cdots,  \text{ with } \la_{\de,n}\ra +\iy \text{ as } n\ra\iy,  \text{ and } \la_{\de,1}\geq C_P^{-1},
	\end{equation}
	as well as 
	\begin{equation}\label{10.19.5}
		0<\la_{\de,1}^{(i)}\leq \la_{\de,2}^{(i)}\leq \la_{\de,3}^{(i)}\leq \cdots,  \text{ with } \la_{\de,n}^{(i)}\ra +\iy \text{ as } n\ra\iy, \text{ and } \la_{\de,1}^{(i)}\geq \f{(1-\al_i)^2}{4}. 
	\end{equation}
	We denote by $\{\Phi_{\de,n}\}_{n\in\N^*}$ and $\{\Phi_{\de,n}^{(i)}\}_{n\in\N^*}$ the corresponding orthonormal bases of $L^2(\Om_\de)$ and $L^2(\La_{\de,i})$, respectively. 
	
	According to the classical Sturm-Liouville oscillation theory and the nodal properties of one-dimensional elliptic operators, for each fixed $\de \in [0, \f{1}{4}]$, the $n$-th eigenfunction $\Phi_{\de,n}^{(i)}$ possesses exactly $n-1$ internal zeros in the interval $\La_{\de,i}$, thereby dividing the domain into exactly $n$ nodal domains with alternating signs. By virtue of the arbitrariness of the complex phase factor, we can assume without loss of generality that $\Phi_{\de,n}^{(i)}$ is real-valued and strictly positive on its first nodal domain adjacent to the left boundary $x_i = \de$.

	Furthermore, the parameterized fractional power spaces $D(\mcA_\de^\theta)\ (\theta \in [0,+\iy))$ and $D((\mcA_\de^{(i)})^{\theta})$ ($i=1,2$) can be defined analogously to Definition \ref{10.18.D1}. The tensor product structure also carries over, yielding the Hilbertian factorization
	\begin{equation*}
		H_0^1(\Om_\de;A_\de) = H_0^1(\La_{\de,1};w_{\de,1})\wh \ots H_0^1(\La_{\de,2};w_{\de,2}).
	\end{equation*}
	Consequently, parallel to Lemma \ref{10.18.L6}, the multi-dimensional eigenfunctions admit the product representation $\Phi_{\de,nm} = \Phi_{\de,n}^{(1)}\Phi_{\de,m}^{(2)}$ for all $n,m\in\N^*$, and every eigenvalue $\la_\de$ of $\mcA_\de$ is uniquely represented as $\la_\de = \la_{\de,nm} = \la_{\de,n}^{(1)} + \la_{\de,m}^{(2)}$.
\end{remark}

\subsection{Weak solution for the shape design}\label{S3.2}

To analyze the impact of boundary perturbations and domain variations, we introduce the following parameterized non-homogeneous Schr\"odinger equation on the localized domain $\Om_\de$ for a fixed shape parameter $\de \in [0, \f{1}{4}]$:
\begin{equation}\label{10.19.6}
	\begin{cases}
		\mathrm{i}\partial_t\varphi_\de - \operatorname{div}(A_\de\nabla\varphi_\de) = f_\de,
		&\text{in } Q_\de,\\
		\varphi_\de = 0,
		&\text{on } \Sigma_\de,\\
		\varphi_\de(0) = \varphi_\de^0,
		&\text{in } \Omega_\de,
	\end{cases}
\end{equation}
where $\vp_\de^0 \in H_0^1(\Om_\de;A_\de)$ and $f_\de \in L^2(Q_\de)$.

\begin{definition}\label{10.19.D2}
	Let $\vp_\de^0\in H_0^1(\Om_\de;A_\de)$ and $f_\de\in L^2(Q_\de)$. We say that a function
	\begin{equation*}
		\vp_\de\in L^2(0,T; H_0^1(\Om_\de;A_\de)), \text{ with }   \pt_t\vp_\de\in L^2(0,T; H^{-1}(\Om_\de;A_\de)) 
	\end{equation*}
	is a weak solution to the equation \eqref{10.19.6} with respect to $(\vp_\de^0,f_\de)$ if:
	
	(1) for each $v\in H_0^1(\Om_\de;A_\de)$ and a.e.~$t\in [0,T]$, 
	\begin{equation*}
		\mathrm{i} \lg \pt_t\vp_\de, v \rg_{H^{-1}(\Om_\de;A_\de), H_0^1(\Om_\de;A_\de)} + a_\de(\vp_\de, v) = (f_\de,v)_{L^2(\Om_\de)},
	\end{equation*}
	
	(2) the initial condition holds in the trace sense, i.e., $\vp_\de(0)=\vp_\de^0$. 
\end{definition}

\begin{theorem}\label{10.19.T1}
	Let $\vp_\de^0\in H_0^1(\Om_\de;A_\de)$ and $f_\de\in L^2(0,T; H_0^1(\Om_\de;A_\de))$. Then, there exists a unique weak solution 
	\begin{equation*}
		\vp_\de\in C([0,T]; H_0^1(\Om_\de;A_\de)), \text{ with }   \pt_t\vp_\de\in L^2(0,T; H^{-1}(\Om_\de;A_\de))
	\end{equation*}
	to the equation \eqref{10.19.6} with respect to $(\vp_\de^0,f_\de)$. Moreover, the solution satisfies the uniform energy estimate:
	\begin{equation}\label{10.19.7}
		\begin{split} 
			&\max_{t\in [0,T]}\|\vp_\de(t)\|_{H_0^1(\Om_\de;A_\de)}+\|\pt_t\vp_\de\|_{L^2(0,T; H^{-1}(\Om_\de;A_\de))}\\
			&\leq C\left(\|\vp_\de^0\|_{H_0^1(\Om_\de;A_\de)}+\|f_\de\|_{L^2(0,T; H_0^1(\Om_\de;A_\de))}\right), 
		\end{split} 
	\end{equation}
	where the positive constant $C$ depends only on $\al_1,\al_2$, and $T$, but is completely independent of $\de \in [0, \f{1}{4}]$. 
\end{theorem}

The proof of Theorem \ref{10.19.T1} is postponed to Appendix \ref{S5}. 

\begin{remark}\label{10.19.R5}
	Let $\vp_\de^0\in H_0^1(\Om_\de;A_\de)$ and $f_\de=0$. Suppose the initial datum $\vp_\de^0$ is expanded in terms of the tensor product basis as
	\begin{equation*}
		\vp_\de^0=\sum_{n,m=1}^\iy \vp_{\de,nm}^0\Phi_{\de,n}^{(1)}\Phi_{\de,m}^{(2)}, \text{ which implies }   \sum_{n,m=1}^\iy \la_{\de, nm}|\vp_{\de, nm}^0|^2<+\iy, 
	\end{equation*}
	where $\vp_{\de, nm}^0=(\vp_\de^0, \Phi_{\de, nm})_{L^2(\Om_\de)}$ for all $n,m\in\N^*$. Then, the weak solution $\vp_\de$ admits the standard Fourier spectral expansion:
	\begin{equation}\label{10.19.8}
		\vp_\de(x_1,x_2,t)=\sum_{n,m=1}^\iy e^{\mathrm{i} \la_{\de, nm}t}\vp_{\de,nm}^0 \Phi_{\de,n}^{(1)}(x_1)\Phi_{\de,m}^{(2)}(x_2). 
	\end{equation}
\end{remark}

\begin{lemma}\label{10.19.L1}
	Let $\de\in (0,\f{1}{4}]$. Assume $\vp_\de^0\in C_0^\iy(\Om_\de)\s D(\mcA_\de)$ and $f_\de=0$. Then the weak solution to the equation \eqref{10.19.6} with respect to $(\vp_\de^0,0)$ satisfies
	\begin{equation}\label{10.19.9}
		\begin{split}
			\vp_\de\in L^2(0,T; D(\mcA_\de)\cap H^2(\Om_\de)), \mbox{ and } \pt_t\vp_\de\in L^2(Q_\de). 
		\end{split}
	\end{equation}
\end{lemma}

\begin{proof}
	From Theorem \ref{10.19.T1}, we obtain $\vp_\de\in L^2(0,T; H_0^1(\Om_\de;A_\de))$. 
	
	(1) We show $\vp_\de\in L^2(0,T; H^2(\Om_\de))$. 
	
	Indeed, from 
	\begin{equation*}
		\begin{split} 
		\pt_{x_ix_i}\Phi_{\de,n}^{(i)}
		&=x_i^{-\al_i}\pt_{x_i}\left(x_i^{\al_i}\pt_{x_i}\Phi_{\de,n}^{(i)}\right)-\al_ix_i^{-1} \pt_{x_i}\Phi_{\de,n}^{(i)}=-x_i^{-\al_i}\mcA_\de^{(i)} \Phi_{\de,n}^{(i)}-\al_ix_i^{-1-\f{\al_i}{2}}\left(x_i^\f{\al_i}{2}\pt_{x_i}\Phi_{\de,n}^{(i)}\right)\\
		&=-x_i^{-\al_i}\la_{\de,n}^{(i)}\Phi_{\de,n}^{(i)}-\al_ix_i^{-1-\f{\al_i}{2}}\left(x_i^\f{\al_i}{2}\pt_{x_i}\Phi_{\de,n}^{(i)}\right), 
		\end{split} 
	\end{equation*}
	we compute
	\begin{equation*}
		\begin{split}
			\left\|\pt_{x_1x_1} \vp_\de(t)\right\|_{L^2(\Om_\de)}^2
			&=\left\|\sum_{n,m=1}^\iy e^{\im \la_{\de,nm}t}\vp_{\de,nm}^0 \left(\pt_{x_1x_1}\Phi_{\de,n}^{(1)}\right)\Phi_{\de,m}^{(2)}\right\|_{L^2(\Om_\de)}^2\\
			&=\left\|\sum_{n,m=1}^\iy e^{\im \la_{\de,nm}t}\vp_{\de,nm}^0 \left(-x_1^{-\al_1}\la_{\de,n}^{(1)}\Phi_{\de,n}^{(1)}-\al_1x_1^{-1-\f{\al_1}{2}}\left(x_1^\f{\al_1}{2}\pt_{x_1}\Phi_{\de,n}^{(1)}\right)\right)\Phi_{\de,m}^{(2)}\right\|_{L^2(\Om_\de)}^2\\
			&\leq 2\de^{-2\al_1}\left\|\sum_{n,m=1}^\iy e^{\im \la_{\de,nm}t}\vp_{\de,nm}^0\la_{\de,n}^{(1)}\Phi_{\de,n}^{(1)}\Phi_{\de,m}^{(2)} \right\|_{L^2(\Om_\de)}^2\\
			&\hspace{4.5mm}+2\al_1^2\de^{-2-\al_1} \left\|\sum_{n,m=1}^\iy e^{\im \la_{\de,nm}t}\vp_{\de,nm}^0\left(x_1^\f{\al_1}{2} \pt_{x_1}\Phi_{\de,n}^{(1)}\right)\Phi_{\de,m}^{(2)}\right\|_{L^2(\Om_\de)}^2\\
			&\leq C\sum_{n,m=1}^\iy \left((\la_{\de,n}^{(1)})^2+\la_{\de,n}^{(1)}\right) |\vp_{\de,nm}^0|^2 
		\end{split}
	\end{equation*}
	by Remark \ref{10.19.R5}, 
	where the positive constant $C$ depends only on $\al_1, \de$. Then, integrating the above inequality over $t \in (0,T)$ yields
	\begin{equation*}
		\left\|\pt_{x_1x_1}\vp_\de\right\|_{L^2(Q_\de)}^2\leq CT\sum_{n,m=1}^\iy \left((\la_{\de,n}^{(1)})^2+\la_{\de,n}^{(1)}\right) |\vp_{\de,nm}^0|^2,
	\end{equation*}
	where the positive constant $C$ depends only on $\al_1$ and $\de$. 
	
	Similar computation we obtain 
	\begin{equation*}
		\left\|\pt_{x_2x_2}\vp_\de\right\|_{L^2(Q_\de)}^2\leq CT\sum_{n,m=1}^\iy \left((\la_{\de,m}^{(2)})^2+\la_{\de,m}^{(2)}\right) |\vp_{\de,nm}^0|^2,
	\end{equation*}
	where the positive constant $C$ depends only on $\al_2$ and $\de$. 
	
	Finally, we compute the mixed second-order spatial derivative. By taking advantage of the separability of the tensor product coordinates and the characterization of the one-dimensional energy norms, we obtain
	\begin{equation*}
		\begin{split}
			\|\pt_{x_1x_2}\vp_\de(t)\|_{L^2(\Om_\de)}^2
			&=\left\|\sum_{n,m=1}^\iy e^{\mathrm{i}\la_{\de,nm}t}\vp_{\de,nm}^0 \left(\pt_{x_1}\Phi_{\de,n}^{(1)}\right)\left(\pt_{x_2}\Phi_{\de,m}^{(2)}\right)\right\|_{L^2(\Om_\de)}^2\\
			&\leq \de^{-\al_1-\al_2}\left\|\sum_{n,m=1}^\iy e^{\mathrm{i} \la_{\de,nm}t} \vp_{\de,nm}^0 \left(x_1^\f{\al_1}{2} \pt_{x_1}\Phi_{\de,n}^{(1)}\right)\left(x_2^{\f{\al_2}{2}}\pt_{x_2}\Phi_{\de,m}^{(2)}\right)\right\|_{L^2(\Om_\de)}^2\\
			&\leq \de^{-\al_1-\al_2} \sum_{n,m=1}^\iy \la_{\de,n}^{(1)}\la_{\de,m}^{(2)}|\vp_{\de,nm}^0|^2 \\
			&\leq \f{1}{4}\de^{-\al_1-\al_2}\sum_{n,m=1}^\iy \left( \la_{\de,n}^{(1)} + \la_{\de,m}^{(2)} \right)^2 |\vp_{\de,nm}^0|^2 \\
			&= \f{1}{4}\de^{-\al_1-\al_2}\sum_{n,m=1}^\iy \la_{\de,nm}^2|\vp_{\de,nm}^0|^2.
		\end{split}
	\end{equation*}
	Integrating this bound over the time interval $(0,T)$ yields
	\begin{equation*}
		\|\pt_{x_1x_2}\vp_\de\|_{L^2(Q_\de)}^2\leq CT\sum_{n,m=1}^\iy \la_{\de,nm}^2|\vp_{\de,nm}^0|^2, 
	\end{equation*}
	where the positive constant $C$ depends only on $\al_1, \al_2$, and $\de$.
	
	Overall, we obtain $\vp_\de\in L^2(0,T; H^2(\Om_\de))$.

	(2) We show $\pt_t\vp_\de\in L^2(Q_\de)$ and $\vp_\de\in L^2(0,T; D(\mcA_\de))$. 
	
	Indeed, we compute
	\begin{equation*}
		\begin{split}
			\|\pt_t\vp_\de(t)\|_{L^2(\Om_\de)}^2
			&=\left\|\sum_{n,m=1}^\iy(\im\la_{\de,nm})e^{\im \la_{\de,nm}t}\vp_{\de,nm}^0 \Phi_{\de,n}^{(1)}\Phi_{\de,m}^{(2)}\right\|_{L^2(\Om_\de)}^2\leq \sum_{n,m=1}^\iy \la_{\de,nm}^2|\vp_{\de,nm}^0|^2, 
		\end{split}
	\end{equation*}
	and then 
	\begin{equation*}
		\|\pt_t\vp_\de\|_{L^2(Q_\de)}^2\leq T\sum_{n,m=1}^\iy \la_{\de,nm}^2|\vp_{\de,nm}^0|^2. 
	\end{equation*}
	Hence $\pt_t\vp_\de\in L^2(Q_\de)$. Same result holds for $\vp_\de\in L^2(0,T; D(\mcA_\de))$. We complete the proof of this lemma. 
\end{proof}

\begin{remark}\label{10.19.R6}
	From the proof of Lemma \ref{10.19.L1}, we observe that the $H^2(\Om_\de)$-regularity of $\varphi_\de$ relies heavily on the fact that the geometric degeneracy is removed by the shape design (i.e., $A_\de$ is uniformly elliptic in $\Om_\de$ for a fixed $\de > 0$). However, for the original degenerate problem \eqref{10.18.24} (corresponding to $\de = 0$), a weak solution $\varphi$ does not belong to $L^2(0,T; H^2(\Om))$ in general. Indeed, the strong degeneracy of the coefficients $x_1^{\alpha_1}$ and $x_2^{\alpha_2}$ prevents the classical elliptic $H^2$-regularity estimates from being applied up to the degenerate boundaries.
\end{remark}

\subsection{Spectral approximation of $\mcA^{(i)}$}\label{S3.3}

Let $i=1,2$ and $\de\in [0,\f{1}{4}]$. We define the parameterized resolvent approximation operator $r_\de^{(i)}$ by
\begin{equation}\label{10.19.11}
	r_\de^{(i)}: L^2(\La_i)\ra L^2(\La_i), \quad f\mapsto r_\de^{(i)}(f)=E_\de^0\left(  u_\de^{(i)}[f]\right), 
\end{equation}
where $E_\de^0$ is the zero extension operator introduced in Definition \ref{10.19.D1}, and $u_\de^{(i)}[f]$ is the unique weak solution to the following boundary value problem on the localized domain:
\begin{equation}\label{10.19.10}
	\begin{cases}
		\mcA_\de^{(i)}\left(u_\de^{(i)}[f]\right)=f_\de, &\mbox{in } \La_{\de,i}, \\
		u_\de^{(i)}[f]=0, &\mbox{on }\pt \La_{\de,i}, 
	\end{cases}
\end{equation}
with $f_\de = f|_{\La_{\de,i}}$. That is, the solution $u_{\de}^{(i)}[f]\in D(\mcA_\de^{(i)})$ satisfies the weak formulation
\begin{equation*}
	\int_{\La_{\de,i}} x_i^{\al_i}\left(\pt_{x_i}u_\de^{(i)}[f]\right)(\pt_{x_i}\ol v)\df x_i=\int_{\La_{\de,i}}f_\de \ol v\df x_i,  \mbox{ for all } v\in H_0^1(\La_{\de,i};w_{\de,i}). 
\end{equation*}
By virtue of the  Poincar\'e  inequality \eqref{10.18.20} and the isometric properties of the zero extension, we can establish the uniform boundedness of $r_\de^{(i)}$ as follows:
\begin{equation*}
	\begin{split}
		\|r_\de^{(i)}(f)\|_{L^2(\La_i)}
		&=\left\|E_\de^0 \left(u_\de^{(i)}[f]\right)\right\|_{L^2(\La_i)} \leq C \|E_\de^0 u_\de^{(i)}[f]\|_{H_0^1(\La_i;w_i)}\\
		&= C\|u_\de^{(i)}[f]\|_{H_0^1(\La_{\de,i};w_{\de,i})} \leq C\|f_\de\|_{L^2(\La_{\de,i})} \leq C\|f\|_{L^2(\La_i)},
	\end{split}
\end{equation*}
where the positive constant $C$ depends exclusively on $\al_i$ and remains strictly independent of $\de \in [0, \f{1}{4}]$. 

Moreover, since the operator $r_\de^{(i)}: L^2(\La_i)\ra L^2(\La_i)$ is bounded and maps $L^2(\La_i)$ into the energy space $H_0^1(\La_i;w_i)$, the compactness of the embedding $H_0^1(\La_i; w_i)\hra L^2(\La_i)$ (Lemma \ref{10.18.L4}) guarantees that $r_\de^{(i)}=(\mcA_\de^{(i)})^{-1}$ is a compact operator. Standard variational arguments further imply that $r_\de^{(i)}$ is a nonnegative and self-adjoint operator on $L^2(\La_i)$.

\begin{lemma}\label{10.19.L2}
	Let $f\in L^2(\La_i)$. Then, the parameterized resolvent operator satisfies the strong convergence
	\begin{equation*}
		r_\de^{(i)}(f)\ra r_0^{(i)}(f)  \mbox{ strongly in } H_0^1(\La_i;w_i)  \mbox{ as } \de\ra 0^+. 
	\end{equation*}
\end{lemma}

\begin{proof}
	Testing the system \eqref{10.19.10} with the weak solution $u_\de^{(i)}[f]\in H_0^1(\La_{\de,i};w_{\de,i})$, we obtain
	\begin{equation*}
		\begin{split}
			\|r_\de^{(i)}(f)\|_{H_0^1(\La_i;w_i)}^2
			&=\int_{\La_i} w_i \left|\pt_{x_i}\left[E_\de^0\left(u_\de^{(i)}[f]\right)\right]\right|^2 \df x_i = \int_{\La_{\de,i}} w_{\de,i} \left|\pt_{x_i} \left(u_\de^{(i)}[f]\right)\right|^2 \df x_i \\
			&=\int_{\La_{\de,i}} f_\de \ol{u_\de^{(i)}[f]} \df x_i \leq \|f_\de\|_{L^2(\La_{\de,i})} \|u_\de^{(i)}[f]\|_{L^2(\La_{\de,i})} \\
			&\leq \f{2}{1-\al_i} \|f\|_{L^2(\La_i)} \|u_\de^{(i)}[f]\|_{H_0^1(\La_{\de,i};w_{\de,i})} = \f{2}{1-\al_i} \|f\|_{L^2(\La_i)} \|r_\de^{(i)}(f)\|_{H_0^1(\La_i;w_i)}, 
		\end{split}
	\end{equation*}
	where we applied the  Poincar\'e inequality \eqref{10.18.20}. This immediately yields the uniform estimate
	\begin{equation*}
		\|r_\de^{(i)}(f)\|_{H_0^1(\La_i;w_i)}\leq \f{2}{1-\al_i}\|f\|_{L^2(\La_i)}. 
	\end{equation*}
	
	Consequently, for any sequence $\de_k \to 0^+$, the boundedness of $\{r_{\de_k}^{(i)}(f)\}$ in the Hilbert space $H_0^1(\La_i;w_i)$ ensures the existence of a subsequence (still denoted by $\de_k$) and a limit function $u_0\in H_0^1(\La_i;w_i)$ such that
	\begin{equation}\label{10.19.12}
		r_{\de_k}^{(i)}(f)\ra u_0 \mbox{ weakly in } H_0^1(\La_i;w_i) \text{ and }  r_{\de_k}^{(i)}(f)\ra u_0  \mbox{ strongly in } L^2(\La_i). 
	\end{equation}
	
	Next, we identify the limit $u_0$. For each test function $v\in C_0^\iy (\La_i)$, we can set $\de_0=\dist(\supp v, \{0\})>0$. Then, for all $\de_k < \de_0$, the support of $v$ is strictly contained in $\La_{\de_k,i}$, which allows us to rewrite the weak formulation \eqref{10.19.10} as
	\begin{equation*}
		\int_{\La_i} x_i^{\al_i} \left(\pt_{x_i} E_{\de_k}^0 u_{\de_k}^{(i)}[f]\right) \ol{\pt_{x_i}v} \df x_i = \int_{\La_i} f \ol v \df x_i. 
	\end{equation*}
	Taking the limit as $k \to \iy$ and utilizing the weak convergence in \eqref{10.19.12}, we obtain
	\begin{equation*}
		\int_{\La_i} x_i^{\al_i} (\pt_{x_i}u_0) \ol{\pt_{x_i}v} \df x_i = \int_{\La_i} f \ol v \df x_i, \mbox{ for all } v\in C_0^\iy(\La_i). 
	\end{equation*}
	By density, this implies that $u_0$ is the unique weak solution to the limiting problem with $\de=0$; hence, $u_0=r_0^{(i)}(f)$, and the uniqueness of the limit guarantees that the entire family converges without passing to a subsequence.
	
	Finally, we establish the convergence of the energy norm. A direct calculation shows
	\begin{equation*}
		\begin{split}
			\lim_{\de\ra 0^+} \int_{\La_i} x_i^{\al_i} |\pt_{x_i} r_\de^{(i)}(f)|^2 \df x_i
			&= \lim_{\de\ra 0^+} \int_{\La_i} f \ol{r_\de^{(i)}(f)} \df x_i \\
			&= \int_{\La_i} f \ol{r_0^{(i)}(f)} \df x_i = \int_{\La_i} x_i^{\al_i}|\pt_{x_i}r_0^{(i)}(f)|^2 \df x_i, 
		\end{split}
	\end{equation*}
	where we used the strong convergence of $r_\de^{(i)}(f)$ in $L^2(\La_i)$. Since weak convergence combined with norm convergence implies strong convergence in any Hilbert space, we conclude that $r_\de^{(i)}(f)\ra r_0^{(i)}(f)$ strongly in $H_0^1(\La_i;w_i)$ as $\de\ra 0^+$. The proof is complete. 
\end{proof}

\begin{remark} 
	Lemma \ref{10.19.L2} established the strong resolvent convergence of the parameterized operators, a property intimately linked to the $\Gamma$-convergence of the associated energy functionals defined on the weighted Sobolev spaces.
\end{remark}

\begin{lemma}[{\cite[Corollaries XI.9.3 and XI.9.4]{Dunford}}]\label{10.19.L3}
	Let $H$ be a Hilbert space, and let $r_1, r_2$ be nonnegative, self-adjoint, and compact operators on $H$. Then, for every $n\in \N^*$, we have 
	\begin{equation*} 
		|\xi_n(r_1)-\xi_n(r_2)|\leq \|r_1-r_2\|_{\mcL(H)},
	\end{equation*}
	where $\xi_n(r_i)\ (i=1,2)$ denotes the $n$-th eigenvalue of $r_i$ ordered according to their magnitudes and counting multiplicities, and $\|\cdot\|_{\mcL(H)}$ denotes the standard operator norm on the space of bounded linear operators.
\end{lemma}

\begin{lemma}\label{10.19.L4}
	Let $i=1,2$. For each $k\in\N^*$, the eigenvalue $\la_{\de,k}^{(i)}$ is monotonically increasing with respect to $\de \in [0, \f{1}{4}]$. Consequently, it is non-increasing as $\de\ra 0^+$ and satisfies the limiting property
	\begin{equation*}
		\la_{\de,k}^{(i)}\ra \la_k^{(i)}  \mbox{ as } \de\ra 0^+.
	\end{equation*}
	In particular, we have the uniform lower bound
	\begin{equation*}
		\la_{\de,k}^{(i)}\geq \la_k^{(i)}   \text{ for all } \de \in \left(0, \f{1}{4}\right].
	\end{equation*}
\end{lemma}

\begin{proof}
	We prove this lemma by the following steps. 
	
	{\it Step 1}. It is well-known that 
	\begin{equation}\label{10.19.13}
		\la_{\de,k}^{(i)}=\min_{E_k\s H_0^1(\La_{\de,i};w_{\de,i})\atop \dim E_k=k}\max_{0\neq v\in E_k}R_\de[v], \mbox{ with } R_\de[v]=\f{\int_{\La_{\de,i}}w_{\de,i}|\pt_{x_i}v|^2\df x_i}{\int_{\La_{\de,i}}|v|^2\df x_i}. 
	\end{equation}
	
	{\it Step 2}. We prove $\la_{\de,k}^{(i)}\leq \la_{\ga,k}^{(i)}$ for all $0\leq \de\leq \ga\leq \f{1}{4}$. 
	
	Indeed, let $E_k\s H_0^1(\La_{\ga,i};w_{\ga,i})$, from Definition \ref{10.19.D1} and Remark \ref{10.19.R3}, we obtain $E_k\s H_0^1(\La_{\de,i};w_{\de,i})$. Since $R_\ga[v]=R_\de[v]$ by Definition \ref{10.19.D1} and Remark \ref{10.19.R3}, we obtain 
	\begin{equation*}
		\la_{\de,k}^{(i)}\leq \la_{\ga,k}^{(i)}, \mbox{ for all } 0\leq \de\leq \ga\leq \f{1}{4}. 
	\end{equation*} 
	
	{\it Step 3}. We prove that $\la_{\de,k}^{(i)}$ is continuous at $\de=0$, i.e., $\la_{\de,k}^{(i)}\ra \la_k^{(i)}$ as $\de\ra 0^+$. 
	
	From Lemma \ref{10.19.L3}, if we can establish the uniform operator convergence
	\begin{equation}\label{10.19.14}
		\|r_\de^{(i)}-r_0^{(i)}\|_{\mcL(L^2(\La_i))}\ra 0 \mbox{ as } \de\ra 0^+, 
	\end{equation}
	then by utilizing the domain monotonicity ($\la_{\de,k}^{(i)} \leq \la_{\f{1}{4},k}^{(i)}$ for all $\de \in [0, \f{1}{4}]$), we can estimate the eigenvalue discrepancy as follows:
	\begin{equation*}
		\left|\la_{\de,k}^{(i)}-\la_k^{(i)}\right| = \la_{\de,k}^{(i)}\la_k^{(i)}\left|\f{1}{\la_{\de,k}^{(i)}}-\f{1}{\la_k^{(i)}}\right| \leq \left(\la_{\f{1}{4},k}^{(i)}\right)^2 \|r_\de^{(i)}-r_0^{(i)}\|_{\mcL(L^2(\La_i))} \ra 0  \mbox{ as } \de\ra 0^+. 
	\end{equation*}
	Hence, it suffices to verify the validity of \eqref{10.19.14}. 
	
	We prove \eqref{10.19.14} by contradiction. Suppose that there exist a constant $\e_0>0$ and a sequence of shape parameters $\de_k \to 0^+$ such that 
	\begin{equation*}
		\|r_{\de_k}^{(i)}-r_0^{(i)}\|_{\mcL(L^2(\La_i))}\geq \e_0, \text{ for all } k\in\N^*. 
	\end{equation*}
	By the definition of the operator norm, there exists a sequence $\{f_{\de_k}\}_{k\in\N^*} \s  L^2(\La_i)$ with $\|f_{\de_k}\|_{L^2(\La_i)}\leq 1$ such that 
	\begin{equation}\label{10.19.17}
		\left\|r_{\de_k}^{(i)}(f_{\de_k})-r_0^{(i)}(f_{\de_k})\right\|_{L^2(\La_i)}\geq \f{1}{2}\e_0. 
	\end{equation}
	The boundedness of $\{f_{\de_k}\}_{k\in\N^*}$ implies that there exists a subsequence (still denoted by itself) and a function $f_0\in L^2(\La_i)$ such that $f_{\de_k}\ra f_0$ weakly in $L^2(\La_i)$ as $k\ra\iy$. 
	
	Furthermore, recalling the uniform boundedness from Lemma \ref{10.19.L2}, we have $\|r_{\de_k}^{(i)}(f_{\de_k})\|_{H_0^1(\La_i;w_i)} \leq \f{2}{1-\al_i} \|f_{\de_k}\|_{L^2(\La_i)} \leq \f{2}{1-\al_i}$. Due to the compact embedding $H_0^1(\La_i;w_i)\hra L^2(\La_i)$, we can extract a further subsequence (still denoted by itself) and a limit function $u_0\in H_0^1(\La_i;w_i)$ such that 
	\begin{equation}\label{10.19.15}
		r_{\de_k}^{(i)}(f_{\de_k})\ra u_0   \mbox{ weakly in } H_0^1(\La_i;w_i)  \text{ and }  r_{\de_k}^{(i)}(f_{\de_k})\ra u_0  \mbox{ strongly in } L^2(\La_i). 
	\end{equation}
	By applying an identical compactness argument to the limiting resolvent sequence, we also find a subsequence (still denoted by itself) and a function $\wt u_0\in H_0^1(\La_i;w_i)$ such that 
	\begin{equation}\label{10.19.16}
		r_0^{(i)}(f_{\de_k})\ra \wt u_0  \mbox{ weakly in } H_0^1(\La_i;w_i) \text{ and }  r_0^{(i)}(f_{\de_k})\ra \wt u_0   \mbox{ strongly in } L^2(\La_i). 
	\end{equation}
	
	Now, let $v\in C_0^\iy(\La_i)$ be an arbitrary test function and set $\de_0=\dist(\supp v,\{0\})>0$. For all sufficiently large $k$ such that $\de_k<\de_0$, we have $\supp v\s  \La_{\de_k,i}$. Taking the weak limit in the variational formulations, we obtain
	\begin{equation*}
		\begin{split} 
			\int_{\La_i} x_i^{\al_i} (\pt_{x_i} u_0)\ol{\pt_{x_i}v}\df x_i
			&= \lim_{k\to\iy} \int_{\La_{i}} x_i^{\al_i} \left(\pt_{x_i} r_{\de_k}^{(i)}(f_{\de_k})\right)\ol{\pt_{x_i}v}\df x_i  = \lim_{k\to\iy} \int_{\La_{\de_k,i}} f_{\de_k}\ol v\df x_i \\
			&= \lim_{k\to\iy} \int_{\La_{\de_k,i}} x_i^{\al_i}\left(\pt_{x_i} u_0^{(i)}[f_{\de_k}]\right) \ol{\pt_{x_i} v}\df x_i = \int_{\La_i} x_i^{\al_i} (\pt_{x_i} \wt u_0)\ol{\pt_{x_i} v}\df x_i,
		\end{split} 
	\end{equation*}
	which implies $u_0 = \wt u_0$ by density. Combining this uniqueness with the strong limits in \eqref{10.19.15} and \eqref{10.19.16}, we deduce that
	\begin{equation*}
		\left\|r_{\de_k}^{(i)}(f_{\de_k})-r_0^{(i)}(f_{\de_k})\right\|_{L^2(\La_i)} \leq \left\|r_{\de_k}^{(i)}(f_{\de_k})-u_0\right\|_{L^2(\La_i)} + \left\|r_0^{(i)}(f_{\de_k})-u_0\right\|_{L^2(\La_i)} \ra 0   \mbox{ as } k\ra\iy. 
	\end{equation*}
	This directly contradicts \eqref{10.19.17}, thereby completing the proof of this lemma. 
\end{proof}

\begin{lemma}\label{10.19.L5}
	Let $i=1,2$. For each $k\in\N^*$, we have 
	\begin{equation*}
		E_\de^0\Phi_{\de,k}^{(i)}\ra \Phi_k^{(i)} \mbox{ strongly in } H_0^1(\La_i;w_i)  \mbox{ as }\de\ra 0^+. 
	\end{equation*}
\end{lemma}

\begin{proof}
	Since $\La_{\de,i}=(\de,1)$, we obtain the $k$-th eigenvalue  $\la_{\de,k}^{(i)}$ is simple. 
	
	From Lemma \ref{10.19.L4}, we obtain 
	\begin{equation*}
		\int_{\La_{i}}x_i^{\al_i} \left|\pt_{x_i}E_\de^0\Phi_{\de,k}^{(i)}\right|^2\df x_i=\int_{\La_{\de,i}}x_i^{\al_i} |\pt_{x_i} \Phi_{\de,k}^{(i)}|^2\df x_i =\la_{\de,k}^{(i)}\leq \la_{\f{1}{4},k}^{(i)}, \mbox{ for all } \de\in \left[0,\f{1}{4}\right], 
	\end{equation*}
	then there exists a subsequence $\{E_{\de_n}^0\Phi_{\de_n,k}^{(i)}\}_{n\in\N^*}$ of $\{E_\de^0\Phi_{\de,k}^{(i)}\}_{\de\in (0,\f{1}{4}]}$, such that 
	\begin{equation}\label{10.20.1}
		E_{\de_n}^0\Phi_{\de_n,k}^{(i)} \ra u_0 \mbox{ weakly in } H_0^1(\La_i;w_i) \mbox{ and } E_{\de_n}^0\Phi_{\de_n,k}^{(i)}\ra u_0 \mbox{ strongly in } L^2(\La_i)
	\end{equation}
	as $\de_n\ra 0^+$. For each $v\in C_0^\iy(\La_i)$, denote $\de_0=\dist(\supp v, \{0\})>0$, then $v\in C_0^\iy(\La_{\de,i})$ for all $\de<(0,\de_0)$,  from 
	\begin{equation*}
		\begin{split} 
		\int_{\La_i} x_i^{\al_i}(\pt_{x_i}E_\de^0\Phi_{\de,k}^{(i)})\ol{\pt_{x_i}v}\df x_i
		&=\int_{\La_{\de,i}}x_i^{\al_i} (\pt_{x_i}\Phi_{\de,k}^{(i)})\ol{\pt_{x_i}v}\df x_i\\
		&=\la_{\de,k}^{(i)}\int_{\La_{\de,i}}\Phi_{\de,k}^{(i)}\ol v\df x_i=\la_{\de,k}^{(i)}\int_{\La_i} (E_\de^0\Phi_{\de,k}^{(i)})\ol v\df x_i
		\end{split} 
	\end{equation*}
	and \eqref{10.20.1} and Lemma \ref{10.19.L4} 
	we obtain 
	\begin{equation*}
		\int_{\La_i} x_i^{\al_i} (\pt_{x_i}u_0)\ol{\pt_{x_i}v}\df x_i=\la_k^{(i)} \int_{\La_i} u_0\ol v\df x_i, \mbox{ for all } v\in C_0^\iy(\La_i), 
	\end{equation*}
	and $\|u_0\|_{L^2(\La_i)}=\lim_{n\ra\iy} \|E_{\de_n}^0\Phi_{\de_n,k}^{(i)}\|_{L^2(\La_i)}=1$. 
	Hence, from the $k$-the eigenvalue $\la_{\de,k}^{(i)}$ is simple by Remarks \ref{10.18.R3} and \ref{10.19.R4}, we obtain 
	\begin{equation}\label{10.20.2}
		u_0=\Phi_k^{(i)}. 
	\end{equation}
	
	Now, from Lemma \ref{10.19.L4} we obtain 
	\begin{equation*}
		\begin{split}
			\int_{\La_i}x_i^{\al_i}|\pt_{x_i}E_{\de_n}^0\Phi_{\de_n,k}^{(i)}|^2\df x_i
			&=\la_{\de_n,k}^{(i)}\ra \la_k^{(i)}=\int_{\La_i} x_i^{\al_i} |\pt_{x_i} \Phi_k^{(i)}|^2\df x_i \mbox{ as } \de_n\ra 0. 
		\end{split}
	\end{equation*}
	Combining  this and \eqref{10.20.1} and  \eqref{10.20.2} and the fact: if $u_n\ra u$ weakly in a Hilbert space and $\|u_n\|_H\ra \|u\|$ then $u_n\ra u$ strongly in $H$, we obtain 
	\begin{equation*}
		E_{\de_n}^0\Phi_{\de_n,k}^{(i)}\ra \Phi_k^{(i)} \mbox{ strongly in } H_0^1(\La_i;w_i)  \mbox{ as } \de_n\ra 0^+. 
	\end{equation*}
	This completes the proof of this lemma. 
\end{proof}

\begin{lemma}\label{10.20.L1}
	Let $i=1,2$. For each $k\in\N^*$, we have 
	\begin{equation}\label{10.20.3}
		\pt_{x_i}\Phi_{\de,k}^{(i)}(1)\ra \pt_{x_i} \Phi_k^{(i)}(1) \mbox{ as } \de \ra 0^+, 
	\end{equation}
	and 
	\begin{equation}\label{10.20.11}
		\left|\pt_{x_i} \Phi_{\de,k}^{(i)}(1)\right|^2\leq C\la_{\de,k}^{(i)}  \mbox{ for all } \de\in \left[0,\f{1}{4}\right], 
	\end{equation}
	where the positive constant $C$ depends only on $\al_i$. 
\end{lemma}

\begin{proof}
	Since 
	\begin{equation*}
		\begin{split} 
		\pt_{x_ix_i} \Phi_{\de,k}^{(i)}
		&=x_i^{-\al_i}\pt_{x_i}\left(x_i^{\al_i}\pt_{x_i}\Phi_{\de,k}^{(i)}\right) -\al_i x_i^{-1} \pt_{x_i}\Phi_{\de,k}^{(i)}=-x_i^{-\al_i} \la_{\de,k}^{(i)}\Phi_{\de,k}^{(i)}-\al_ix_i^{-1-\f{\al_i}{2}}\left(x_i^\f{\al_i}{2}\pt_{x_i}\Phi_{\de,k}^{(i)}\right), 
		\end{split} 
	\end{equation*}
	and $\Phi_{\de,k}^{(i)}, x_i^\f{\al_i}{2}\pt_{x_i}\Phi_{\de,k}^{(i)}\in L^2(\La_{\de,i})$, 
	we obtain $\Phi_{\de,k}^{(i)}\in H^2(\La_{\f{1}{4},i})$. 
	
	Let $\zeta=\zeta(x_i)\in C^\iy(\R), 0\leq \zeta\leq 1$ such that 
	\begin{equation*}
		\begin{split}
			\zeta=0 \mbox{ on }\left(-\iy, \f{1}{4}\right),\quad \zeta=1 \mbox{ on }\left(\f{1}{2},+\iy\right), \quad |\zeta'|+|\zeta''|\leq C, 
		\end{split}
	\end{equation*}
	where the positive constant $C$ is absolute. 
	
	Testing $\mcA_\de^{(i)}\Phi_{\de,k}^{(i)}=\la_{\de,k}^{(i)}\Phi_{\de,k}^{(i)}$ by $\zeta \pt_{x_i} \Phi_{\de,k}^{(i)}$, integration by parts on $\La_{\de,i}$, since $\supp (\zeta \Phi_{\de,k}^{(i)})\s [\f{1}{4},1]$ and $\Phi_{\de,k}^{(i)}\in L^2(\La_{\f{1}{4},i})$, we obtain 
	\begin{equation}\label{10.20.4}
		\begin{split}
			\f{1}{2}\left|\pt_{x_i}\Phi_{\de,k}^{(i)}\right|^2(1)
			&=\f{1}{2}\int_{\La_{\de,i}}(\pt_{x_i}\zeta)x_i^{\al_i}|\pt_{x_i}\Phi_{\de,k}^{(i)}|^2\df x_i-\f{\al_i}{2}\int_{\La_{\de,i}}\zeta x_i^{\al_i-1}|\pt_{x_i}\Phi_{\de,k}^{(i)}|^2\df x_i\\
			&\hspace{4.5mm}+\f{1}{2}\la_{\de,k}^{(i)}\int_{\La_{\de,i}}(\pt_{x_i}\zeta)|\Phi_{\de,k}^{(i)}|^2\df x_i. 
		\end{split}
	\end{equation}
	Together with this and $\int_{\La_{\de,i}}x_i^{\al_i}|\pt_{x_i}\Phi_{\de,k}^{(i)}|^2\df x_i=\la_{\de,k}^{(i)}$, we get 
	\begin{equation*}
		\begin{split}
			\left|\pt_{x_i}\Phi_{\de,k}^{(i)}\right|^2(1)\leq C \la_{\de,k}^{(i)}, 
		\end{split}
	\end{equation*}
	where the positive constant $C$ depends only on $\al_i$. Moreover, from \eqref{10.20.4} and Lemmas \ref{10.19.L4} and  \ref{10.19.L5}, we obtain 
	\begin{equation*}
		\pt_{x_i}\Phi_{\de,k}^{(i)}(1)\ra \pt_{x_i}\Phi_k^{(i)}(1) \mbox{ as } \de\ra 0^+. 
	\end{equation*}
	We complete the proof. 
\end{proof}

\subsection{Approximation}\label{S3.4}

Let $\vp^0\in C_0^\iy(\Om)$.  Denote 
\begin{equation}\label{10.20.6}
	\de\in [0,\de_0], \mbox{ with }\de_0=\f{1}{2}\min\left\{\f{1}{4}, \dist(\supp \vp^0, \pt\Om)\right\}. 
\end{equation}
Then $\vp^0\in C_0^\iy(\Om_\de)$ for all $\de\in [0,\de_0]$. 

For $\de\in [0,\de_0]$, we have the spectral expansions
\begin{equation}\label{10.20.8}
	\begin{split}
		\vp^0=\sum_{n,m=1}^\iy \vp_{\de, nm}^0\Phi_{\de,nm}, \mbox{ with } \vp_{\de,nm}^0=(\vp^0, \Phi_{\de,nm})_{L^2(\Om_\de)}. 
	\end{split}
\end{equation}
and  the unique weak solution to the perturbed problem \eqref{10.19.6} corresponding to the initial datum $(\vp^0, 0)$ admits the explicit spectral representation
\begin{equation}\label{10.20.7}
	\vp_\de(t) = \sum_{n,m=1}^\iy e^{\mathrm{i}\la_{\de, nm}t} \vp_{\de, nm}^0 \Phi_{\de,nm} = \sum_{n,m=1}^\iy e^{\mathrm{i} \la_{\de,nm}t} \vp_{\de,nm}^0 \Phi_{\de,n}^{(1)} \Phi_{\de,m}^{(2)}.
\end{equation}
This formulation naturally aligns with the parallel spectral expansions previously established in \eqref{10.19.2} and \eqref{10.19.8}. Moreover, $\vp_0=\vp$.

We now define a uniform higher-order regularity bound $M$ across the parameterized family of subdomains. Let
\begin{equation}\label{10.20.9}
	M = \sup_{\de\in [0,\de_0]} M_\de+1,   \text{ with }  M_\de = \sum_{j=0}^4 \sum_{n,m=1}^\iy \la_{\de,nm}^j |\vp_{\de,nm}^0|^2.
\end{equation}
By virtue of Remark \ref{10.19.R4} and the geometric threshold \eqref{10.20.6}, the initial datum $\vp^0$ possesses compact support strictly contained within the localized domain $\Om_\de$. Consequently, the boundary singularity is avoided, and the standard elliptic regularity elliptic estimates can be applied uniformly. For instance, the second-order relation satisfies
\begin{equation*}
	\|\mcA_\de\vp^0\|_{L^2(\Om_\de)}^2 = \int_{\Om_\de} \left|\operatorname{div}(A_\de\nabla \vp^0)\right|^2 \df x \leq C \sum_{i=0}^2 \|D^i\vp^0\|_{L^\iy(\Om)}^2, 
\end{equation*}
where the positive constant $C$ depends exclusively on $\al_1, \al_2$, and $\de_0$, but remains completely independent of the domain perturbation parameter $\de\in [0,\de_0]$.
More generally, by expressing the spectral components in terms of the fractional powers of the operator $\mcA_\de$, we obtain the following uniform upper bound:
\begin{equation*}
	M_\de = \sum_{j=0}^4 \|\mcA_\de^{\f{j}{2}}\vp^0\|_{L^2(\Om_\de)}^2 \leq C \sum_{j=0}^4 \|D^j\vp^0\|_{L^\iy(\Om)}^2,
\end{equation*}
where the positive constant $C$ depends exclusively on $\al_1, \al_2$, and $\de_0$, but remains completely independent of the domain perturbation parameter $\de\in [0,\de_0]$.

\begin{lemma}\label{10.20.L2}
	Let $\vp^0\in C_0^\iy(\Om)$ and $\de\in [0,\de_0]$, where $\de_0$ is the geometric threshold parameter defined in \eqref{10.20.6}. Let $\vp_\de$ be the unique weak solution to the perturbed problem \eqref{10.19.6} corresponding to the data $(\vp^0,0)$. Then, the conormal derivative of the solution satisfies the boundary trace regularity
	\begin{equation*}
		\f{\pt\vp_\de}{\pt \nu_{A_\de}}\in L^2\left(0,T; L^2(\Ga_{\de,1}^{1}\cup \Ga_{\de,2}^{1})\right). 
	\end{equation*}
\end{lemma}

\begin{proof}
	It suffices to establish the boundary regularity on the vertical component of the non-degenerate boundary, namely, $\f{\pt\vp_\de}{\pt\nu_{A_\de}}\in L^2(0,T; L^2(\Ga_{\de,1}^{1}))$. An analogous proof directly establishes the parallel assertion for the horizontal segment $\Ga_{\de,2}^{1}$. 
	
	Since $\f{\pt\vp_\de}{\pt\nu_{A_\de}}=\nabla \vp_\de\cdot A_\de\nu=\pt_{x_1}\vp_\de(1,\cdot,\cdot)$ on $\Ga_{\de,1}^1\ts (0,T)$, and noting that $\la_{\de,nm}=\la_{\de,n}^{(1)}+\la_{\de,m}^{(2)}\geq \la_{\de,n}^{(1)}$ for all $n,m\in\N^*$, we differentiate the spectral expansion \eqref{10.20.7} with respect to $x_1$ at the boundary. By exploiting the orthonormality of $\{\Phi_{\de,m}^{(2)}\}$ in $L^2(\Ga_{\de,1}^1)$ and applying the  Cauchy-Schwarz inequality, we compute
	\begin{equation*}
		\begin{split}
			\left\|\pt_{x_1}\vp_\de(1,\cdot,t)\right\|_{L^2(\Ga_{\de,1}^1)}^2
			&=\left\|\sum_{n,m=1}^\iy e^{\mathrm{i}\la_{\de,nm}t}\vp_{\de, nm}^0 \left(\pt_{x_1}\Phi_{\de,n}^{(1)}(1)\right)\Phi_{\de,m}^{(2)} \right\|_{L^2(\Ga_{\de,1}^1)}^2\\
			&=\sum_{m=1}^\iy \left|\sum_{n=1}^\iy e^{\mathrm{i} \la_{\de,nm}t}\vp_{\de,nm}^0 \pt_{x_1}\Phi_{\de,n}^{(1)}(1)\right|^2\\
			&=\sum_{m=1}^\iy \left|\sum_{n=1}^\iy e^{\mathrm{i}\la_{\de,nm}t}\vp_{\de,nm}^0\la_{\de,nm}^\f{3}{2} \left[\la_{\de,nm}^{-\f{3}{2}} \pt_{x_1}\Phi_{\de,n}^{(1)}(1)\right]\right|^2\\
			&\leq \sum_{m=1}^\iy \left(\sum_{n=1}^\iy \la_{\de,nm}^3|\vp_{\de,nm}^0|^2\right)\left(\sum_{n=1}^\iy \la_{\de,nm}^{-3} \left|\pt_{x_1}\Phi_{\de,n}^{(1)}(1)\right|^2\right).
		\end{split}
	\end{equation*}
	Invoking the uniform spectral boundary trace estimate from Lemma \ref{10.20.L1}, we have $\left|\pt_{x_i}\Phi_{\de,k}^{(i)}(1)\right|^2 \leq C \la_{\de,k}^{(i)} \leq C \la_{\de,nm}$, where the positive constant $C$ depends only on $\al_i$. Substituting this bound into the second series yields
	\begin{equation*}
		\begin{split}
			\left\|\pt_{x_1}\vp_\de(1,\cdot,t)\right\|_{L^2(\Ga_{\de,1}^1)}^2
			&\leq C \sum_{m=1}^\iy \left(\sum_{n=1}^\iy \la_{\de,nm}^3|\vp_{\de,nm}^0|^2\right)\left(\sum_{n=1}^\iy \la_{\de,nm}^{-2}\right)\\
			&\leq C\left(\sum_{n=1}^\iy (\la_{\de,n}^{(1)})^{-2}\right) \sum_{n,m=1}^\iy \la_{\de,nm}^3|\vp_{\de,nm}^0|^2 < +\iy, 
		\end{split}
	\end{equation*}
	where we used the fact that the spectral Dirichlet series $\sum_{n=1}^\iy (\la_{\de,n}^{(1)})^{-2}$ converges uniformly due to the standard growth order of eigenvalues by Remark \ref{10.18.R3} and Lemma \ref{10.19.L4}. 
	
	Then, integrating this pointwise bound over the time interval $(0,T)$ and recalling the definition of the higher-order spectral bound $M$ in \eqref{10.20.9}, we deduce that
	\begin{equation*}
		\left\|\pt_{x_1}\vp_\de(1,\cdot,\cdot)\right\|_{L^2(0,T;L^2(\Ga_{\de,1}^1))}^2 \leq CT\left(\sum_{n=1}^\iy (\la_{\de,n}^{(1)})^{-2}\right) \sum_{n,m=1}^\iy \la_{\de,nm}^3|\vp_{\de,nm}^0|^2 \leq C T M. 
	\end{equation*}
	This confirms that $\f{\pt\vp_\de}{\pt\nu_{A_\de}}\in L^2(0,T; L^2(\Ga_{\de,1}^{1}))$, and the proof is complete. 
\end{proof}

\begin{lemma}\label{10.20.L3}
	Let $\vp^0\in C_0^\iy(\Om)$ and $\de\in [0,\de_0]$, where $\de_0$ is the geometric threshold parameter defined in \eqref{10.20.6}. Let $\vp_\de$ be the unique weak solution to the perturbed problem \eqref{10.19.6} corresponding to the data $(\vp^0,0)$. Then, the solution satisfies the $L^2$ and $H_0^1$ energy conservation laws:
	\begin{equation}\label{10.20.12}
		\|\vp_\de(t)\|_{L^2(\Om_\de)} = \|\vp^0\|_{L^2(\Om_\de)}  \text{ and }  \|\vp_\de(t)\|_{H_0^1(\Om_\de;A_\de)} = \|\vp^0\|_{H_0^1(\Om_\de;A_\de)}
	\end{equation}
	for all $t\in [0,T]$. 
\end{lemma}

\begin{proof}
	By invoking the spectral expansions \eqref{10.20.8} and \eqref{10.20.7}, the orthonormality of the multi-dimensional tensor product basis $\{\Phi_{\de,nm}\}_{n,m\in\N^*}$ combined with the property $|e^{\mathrm{i}\la_{\de,nm}t}| = 1$ allows us to apply Parseval's identity pointwise in time. Specifically, for the $L^2(\Om_\de)$-norm, we compute
	\begin{equation*}
		\|\vp_\de(t)\|_{L^2(\Om_\de)}^2 = \left\|\sum_{n,m=1}^\iy e^{\mathrm{i} \la_{\de,nm}t}\vp_{\de,nm}^0 \Phi_{\de,nm}\right\|_{L^2(\Om_\de)}^2 = \sum_{n,m=1}^\iy |\vp_{\de,nm}^0|^2 = \|\vp^0\|_{L^2(\Om_\de)}^2.
	\end{equation*}
	Similarly, for the equivalent energy norm on $H_0^1(\Om_\de;A_\de)$ characterized by the sesquilinear form $a_\de(\cdot,\cdot)$, we readily obtain
	\begin{equation*}
		\|\vp_\de(t)\|_{H_0^1(\Om_\de;A_\de)}^2 = \sum_{n,m=1}^\iy \la_{\de,nm}|\vp_{\de,nm}^0|^2 = \|\vp^0\|_{H_0^1(\Om_\de;A_\de)}^2. 
	\end{equation*}
	This completes the proof of Lemma \ref{10.20.L3}.
\end{proof}

\begin{definition}\label{10.20.D_spacetime}
	Let $y\in L^2(Q_\de)$. We define the space-time zero extension operator $E_\de^0$ by
	\begin{equation*}
		E_\de^0 y(x,t)=
		\begin{cases}
			y(x,t), &(x,t)\in Q_\de,\\
			0, &(x,t)\in Q-Q_\de.
		\end{cases}
	\end{equation*}
	Then, $E_\de^0y \in L^2(Q)$ and acts as an isometry, satisfying $\|E_\de^0 y\|_{L^2(Q)}=\|y\|_{L^2(Q_\de)}$. Furthermore, the following structural extensions hold true:
	
	(1) If $y\in L^2(0,T; H_0^1(\Om_\de;A_\de))$, then $E_\de^0 y\in L^2(0,T; H_0^1(\Om;A))$ with the preservation of the energy norm:
	\begin{equation*}
		\|E_\de^0 y\|_{L^2(0,T; H_0^1(\Om;A))}=\|y\|_{L^2(0,T; H_0^1(\Om_\de;A_\de))}.
	\end{equation*}
	
	(2) If $y\in H^1(0,T; L^2(\Om_\de))$, then $E_\de^0 y\in H^1(0,T; L^2(\Om))$. Moreover, the time derivative commutes with the extension operator, i.e., $\pt_t (E_\de^0 y) = E_\de^0 (\pt_t y)$, and the norm satisfies
	\begin{equation*}
		\|E_\de^0 y\|_{H^1(0,T; L^2(\Om))}=\|y\|_{H^1(0,T; L^2(\Om_\de))}.
	\end{equation*}
\end{definition}

\begin{theorem}\label{10.20.T1}
	Let $\vp^0\in C_0^\iy(\Om)$ and $\de\in [0,\de_0]$, where $\de_0$ is the geometric threshold parameter defined in \eqref{10.20.6}. Let $\vp_\de$ be the unique weak solution to the perturbed problem \eqref{10.19.6} corresponding to the data $(\vp^0,0)$. Then, up to a subsequence as $\de\ra 0^+$, 
	\begin{equation}\label{10.20.13}
		\begin{split}
			E_\de^0 \vp_\de\ra \vp\equiv \vp_0 \mbox{ strongly in } L^2(0,T; H_0^1(\Om;A)), 
		\end{split}
	\end{equation}
	and 
	\begin{equation}\label{10.20.14}
		\begin{split}
			\f{\pt E_\de^0\vp_\de}{\pt \nu_{A}}\ra \f{\pt\vp}{\pt \nu_A} \mbox{ strongly in } L^2\left(0,T; L^2(\Ga_1^1\cup \Ga_2^1)\right).  
		\end{split}
	\end{equation}
\end{theorem}

\begin{proof}
	We establish the verification of this theorem by dividing the argument into the following structural steps. 
	
	{\it Step 1}. We first prove the validity of the strong convergence in \eqref{10.20.13}. 
	
	Let $\e>0$ be arbitrary. By invoking the uniform growth of the one-dimensional eigenvalues, we can choose a sufficiently large cutoff index $k\in\N^*$ such that for all $i=1,2$, 
	\begin{equation}\label{10.20.15}
		\la_n^{(i)}\geq \max\left\{1, M\e^{-1}\right\} \text{ for all } n\geq k, \text{ and }  \sum_{n=k+1}^\iy \f{1}{\la_n^{(i)}} < \f{\e}{M}.
	\end{equation}
	By virtue of the domain monotonicity guaranteed by Lemma \ref{10.19.L4}, the perturbed eigenvalues satisfy the parallel bounds uniformly across all shape parameters $\de\in [0,\de_0]$:
	\begin{equation*}
		\la_{\de,n}^{(i)}\geq \max\left\{1, M\e^{-1}\right\} \text{ for all } n\geq k,   \text{ and }   \sum_{n=k+1}^\iy \f{1}{\la_{\de,n}^{(i)}} < \f{\e}{M}.
	\end{equation*}
	To isolate the high-frequency components of the spectrum, we define the index remainder set
	\begin{equation}\label{10.20.16}
		A(k)=\left\{(n,m)\in\N^2\colon n\geq k+1 \mbox{ or } m\geq k+1\right\}. 
	\end{equation}
	A direct application of the elementary algebraic inequality $|a+b+c|^2 \leq 3|a|^2 + 3|b|^2 + 3|c|^2$ yields the following split for the energy norm error:
	\begin{equation*}
		\begin{split}
			&\|E_\de^0\vp_\de(t)-\vp(t)\|_{H_0^1(\Om;A)}^2\\
			&=\left\|\sum_{n,m=1}^\iy e^{\mathrm{i} \la_{\de,nm}t} \vp_{\de,nm}^0 E_\de^0\Phi_{\de,nm} - \sum_{n,m=1}^\iy e^{\mathrm{i}\la_{nm}t}\vp_{nm}^0 \Phi_{nm} \right\|_{H_0^1(\Om;A)}^2\\
			&\leq 3\left\|\sum_{n,m=1}^k e^{\mathrm{i} \la_{\de,nm}t} \vp_{\de,nm}^0 E_\de^0\Phi_{\de,nm}- \sum_{n,m=1}^k e^{\mathrm{i}\la_{nm}t}\vp_{nm}^0\Phi_{nm} \right\|_{H_0^1(\Om;A)}^2\\
			&\quad + 3\left\|\sum_{(n,m)\in A(k)}e^{\mathrm{i} \la_{\de,nm}t}\vp_{\de,nm}^0 E_\de^0\Phi_{\de,nm}\right\|_{H_0^1(\Om;A)}^2 + 3\left\|\sum_{(n,m)\in A(k)}e^{\mathrm{i}\la_{nm}t}\vp_{nm}^0 \Phi_{nm}\right\|_{H_0^1(\Om;A)}^2\\
			&=: B_1 + B_2 + B_3. 
		\end{split}
	\end{equation*}
	
	We now estimate the high-frequency perturbed component $B_2$. For any shape parameter $\de\in [0,\de_0]$, by decomposing the set $A(k)$ and using the relation $\la_{\de,nm}=\la_{\de,n}^{(1)}+\la_{\de,m}^{(2)}$, we have
	\begin{equation*}
		\begin{split}
			B_2
			&= 3 \sum_{(n,m)\in A(k)} \la_{\de,nm}|\vp_{\de,nm}^0|^2 \leq 3 \sum_{n\in\N^*,m\geq k+1} \la_{\de,nm}|\vp_{\de,nm}^0|^2 + 3 \sum_{n\geq k+1, m\in\N^*} \la_{\de,nm}|\vp_{\de,nm}^0|^2 \\
			&\leq 3 \sum_{n\in\N^*, m\geq k+1} \la_{\de,nm}\la_{\de,m}^{(2)}|\vp_{\de,nm}^0|^2 (\la_{\de,m}^{(2)})^{-1} + 3 \sum_{n\geq k+1,m\in\N^*} \la_{\de,nm}\la_{\de,n}^{(1)}|\vp_{\de,nm}^0|^2 (\la_{\de,n}^{(1)})^{-1} \\
			&\leq 3 \f{\e}{M} \sum_{n,m=1}^\iy \la_{\de,nm}^2|\vp_{\de,nm}^0|^2 + 3 \f{\e}{M} \sum_{n,m=1}^\iy \la_{\de,nm}^2|\vp_{\de,nm}^0|^2 \leq 6\e,
		\end{split}
	\end{equation*}
	where we invoked the higher-order spectral uniformity relation \eqref{10.20.9}. 
	
	Now, to control the low-frequency component $B_1$ involving finite terms, we expand the individual spectral differences using a standard algebraic telescopic splitting:
	\begin{equation*}
		\begin{split}
			&e^{\mathrm{i} \la_{\de,nm}t} \vp_{\de,nm}^0 E_\de^0(\Phi_{\de,nm}) - e^{\mathrm{i} \la_{nm}t} \vp_{nm}^0 \Phi_{nm}\\
			&= \left(e^{\mathrm{i} \la_{\de,nm}t} - e^{\mathrm{i}\la_{nm}t}\right) \vp_{\de,nm}^0 E_\de^0\Phi_{\de,nm} + e^{\mathrm{i}\la_{nm}t} \left(\vp_{\de,nm}^0 - \vp_{nm}^0\right) E_\de^0\Phi_{\de,nm} \\
			&\quad + e^{\mathrm{i}\la_{nm}t} \vp_{nm}^0 \left( E_\de^0\Phi_{\de,nm} - \Phi_{nm} \right).
		\end{split}
	\end{equation*}
	
	For the first component, by the mean value theorem, the time-harmonic phase difference is uniformly bounded by the eigenvalue discrepancy for each $t \in [0,T]$:
	\begin{equation}\label{10.20.18}
		\left|e^{\mathrm{i}\la_{\de,nm}t} - e^{\mathrm{i} \la_{nm}t}\right| \leq T |\la_{\de,nm}-\la_{nm}| \ra 0 \text{ as } \de \ra 0^+,
	\end{equation}
	which vanishes according to the eigenvalue continuity in Lemma \ref{10.19.L4}. For the multi-dimensional eigenfunctions, since 
	\begin{equation*}
		E_\de^0\Phi_{\de,nm} - \Phi_{nm} = (E_\de^0\Phi_{\de,n}^{(1)})\left(E_\de^0\Phi_{\de,m}^{(2)}- \Phi_m^{(2)}\right) + \Phi_m^{(2)}\left(E_\de^0\Phi_{\de,n}^{(1)} - \Phi_n^{(1)}\right), 
	\end{equation*} invoking the one-dimensional strong convergence from Lemma \ref{10.19.L5} alongside the uniform boundedness guarantees that $E_\de^0\Phi_{\de,nm} \ra \Phi_{nm}$ strongly in $H_0^1(\Om;A)$ as $\de \ra 0^+$. 
	
	Consequently, for the Fourier coefficients, the complex-valued $L^2$-inner product difference satisfies
	\begin{equation}\label{10.20.19}
		\vp_{\de,nm}^0 - \vp_{nm}^0 = (\vp^0, \Phi_{\de,nm})_{L^2(\Om_\de)} - (\vp^0, \Phi_{nm})_{L^2(\Om)} = \int_\Om \vp^0 \left( \ol{E_\de^0\Phi_{\de,nm} - \Phi_{nm}} \right) \df x \ra 0.
	\end{equation}
	Since the sum in $B_1$ contains only a finite number of terms ($1 \leq n,m \leq k$), the convergence of each individual component implies that there exists a sufficiently small threshold $\de_1 > 0$ such that for all $\de\in (0,\de_1)$, we have $B_1 \leq \e$. 
	
	By selecting $B_3 \leq 6\e$ analogously to $B_2$ via the tail decay of the baseline spectral series, and combining the bounds for $B_1, B_2$, and $B_3$, we arrive at
	\begin{equation*}
		\|E_\de^0\vp_\de(t)-\vp(t)\|_{H_0^1(\Om;A)}^2 \leq B_1 + B_2 + B_3 \leq \e + 6\e + 6\e = 13\e
	\end{equation*}
	whenever $\de\in (0,\de_1)$. By the arbitrariness of $\e > 0$, this completes the proof of the strong convergence in the energy space.

	{\it Step 2}. We next prove the strong convergence of the normal derivatives on the non-degenerate boundaries as specified in \eqref{10.20.14}. It suffices to establish the convergence on the vertical component of the boundary, namely,
	\begin{equation}\label{10.20.17}
		\f{\pt E_\de^0\vp_\de}{\pt\nu_A} \ra \f{\pt\vp}{\pt \nu_A} \mbox{ strongly in } L^2\left(0,T; L^2(\Ga_{1}^{1})\right)  \text{ as } \de \ra 0^+.
	\end{equation}
	Once this is established, a parallel argument immediately yields the identical strong convergence on the horizontal segment of the boundary, i.e., $\f{\pt E_\de^0\vp_\de}{\pt\nu_A} \ra \f{\pt\vp}{\pt \nu_A}$ strongly in $L^2\left(0,T; L^2(\Ga_{2}^{1})\right)$.

	Let $\e>0$ be arbitrary. Let $k\in\N^*$ be the specific spectral cutoff index chosen to satisfy the decay conditions in \eqref{10.20.15}, and let $A(k)$ be the associated high-frequency index remainder set defined in \eqref{10.20.16}. 
	
	Since $\f{\pt (E_\de^0\vp_{\de})}{\pt \nu_A}=\nabla (E_\de^0\vp_\de) \cdot A\nu=[\pt_{x_1}(E_\de^0\vp_\de)](1,\cdot,\cdot)$ holds on $\Ga_1^1\ts (0,T)$, by exploiting the orthonormality of the tensor product components and applying the algebraic inequality $|a+b+c|^2 \leq 3|a|^2 + 3|b|^2 + 3|c|^2$, we can partition the boundary trace error into low- and high-frequency components as follows:
	\begin{equation*}
		\begin{split}
			&\left\|\f{\pt (E_\de^0\vp_\de)}{\pt\nu_A}(t)-\f{\pt\vp}{\pt \nu_A}(t)\right\|_{L^2(\Ga_1^1)}^2\\
			&=\left\|\sum_{n,m=1}^\iy e^{\mathrm{i}\la_{\de,nm}t}\vp_{\de,nm}^0 \pt_{x_1}\Phi_{\de,n}^{(1)}(1) E_\de^0\Phi_{\de,m}^{(2)}-\sum_{n,m=1}^\iy e^{\mathrm{i}\la_{nm}t}\vp_{nm}^0 \pt_{x_1}\Phi_n^{(1)}(1)\Phi_m^{(2)}\right\|_{L^2(\Ga_1^1)}^2\\
			&\leq 3\left\|\sum_{n,m=1}^k e^{\mathrm{i}\la_{\de,nm}t}\vp_{\de,nm}^0 \pt_{x_1}\Phi_{\de,n}^{(1)}(1) E_\de^0\Phi_{\de,m}^{(2)}-\sum_{n,m=1}^k e^{\mathrm{i}\la_{nm}t}\vp_{nm}^0 \pt_{x_1}\Phi_n^{(1)}(1)\Phi_m^{(2)}\right\|_{L^2(\Ga_1^1)}^2\\
			&\quad +3\left\|\sum_{(n,m)\in A(k)} e^{\mathrm{i}\la_{\de,nm}t}\vp_{\de,nm}^0 \pt_{x_1}\Phi_{\de,n}^{(1)}(1) E_\de^0\Phi_{\de,m}^{(2)}\right\|_{L^2(\Ga_1^1)}^2\\
			&\quad +3\left\|\sum_{(n,m)\in A(k)} e^{\mathrm{i}\la_{nm}t}\vp_{nm}^0 \pt_{x_1}\Phi_n^{(1)}(1)\Phi_m^{(2)}\right\|_{L^2(\Ga_1^1)}^2\\
			&=: D_1+D_2+D_3. 
		\end{split}
	\end{equation*}
	
	To bound the high-frequency perturbed trace component $D_2$, we first establish a uniform convergence relation for the weighted spectral weights. Recalling the relation $\la_{\de,nm}=\la_{\de,n}^{(1)}+\la_{\de,m}^{(2)}$ and the growth conditions in \eqref{10.20.15}, we observe that
	\begin{equation*}
		\begin{split}
			\sup_{m\in\N^*}\sum_{n: (n,m)\in A(k)}\la_{\de,nm}^{-3}\la_{\de,n}^{(1)}
			&\leq \sup_{m\in\N^*}\sum_{n: (n,m)\in A(k)}\la_{\de,nm}^{-2}=\sup_{m\in\N^*}\sum_{n: (n,m)\in A(k)}\la_{\de,nm}^{-1}\f{1}{\la_{\de,n}^{(1)}+\la_{\de,m}^{(2)}}\\
			&\leq M^{-1}\e \sup_{m\in\N^*}\sum_{n:(n,m)\in A(k)} \la_{\de,nm}^{-1}\leq M^{-1}\e \sum_{n\in\N^*}\f{1}{\la_n^{(1)}}\leq CM^{-1}\e.
		\end{split}
	\end{equation*}
	
	Next, by applying the discrete Cauchy-Schwarz inequality to the spectral coefficients of $D_2$ and invoking the uniform sublinear boundary estimate from Lemma \ref{10.20.L1}, we deduce that for each $t \in [0,T]$:
	\begin{equation*}
		\begin{split}
			D_2
			&=3\sum_{m=1}^\iy \left|\sum_{n: (n,m)\in A(k)} e^{\mathrm{i}\la_{\de,nm}t}\vp_{\de,nm}^0 \pt_{x_1}\Phi_{\de,n}^{(1)}(1) \right|^2\\
			&=3\sum_{m=1}^\iy \left|\sum_{n: (n,m)\in A(k)} e^{\mathrm{i}\la_{\de,nm}t}\la_{\de,nm}^\f{3}{2}\vp_{\de,nm}^0 \pt_{x_1}\Phi_{\de,n}^{(1)}(1) \la_{\de,nm}^{-\f{3}{2}} \right|^2\\
			&\leq 3\sum_{m=1}^\iy \left(\sum_{n: (n,m)\in A(k)}\la_{\de,nm}^3|\vp_{\de,nm}^0|^2\right)\left(\sum_{n: (n,m)\in A(k)}\la_{\de,nm}^{-3}\left|\pt_{x_1}\Phi_{\de,n}^{(1)}(1)\right|^2\right)\\
			&\leq C\sum_{m=1}^\iy \left(\sum_{n: (n,m)\in A(k)}\la_{\de,nm}^3|\vp_{\de,nm}^0|^2\right)\left(\sum_{n: (n,m)\in A(k)}\la_{\de,nm}^{-3}\la_{\de,n}^{(1)}\right)\leq C\e,
		\end{split}
	\end{equation*}
	where the positive constant $C$ remains completely independent of $\de \in [0, \de_0]$.

	Parallelly, by exploiting the tail decay of the baseline spectral series corresponding to $\de=0$, an identical argument yields
	\begin{equation*}
		D_3 \leq C\e,
	\end{equation*}
	where the positive constant $C$ is independent of the cutoff index $k$.
	
	Finally, to control the low-frequency boundary trace component $D_1$, which contains only a finite number of terms ($1 \leq n,m \leq k$), we implement a standard algebraic telescopic splitting for the individual spectral terms:
	\begin{equation*}
		\begin{split}
			&e^{\mathrm{i}\la_{\de,nm}t}\vp_{\de,nm}^0 \pt_{x_1}\Phi_{\de,n}^{(1)}(1) E_\de^0(\Phi_{\de,m}^{(2)}) - e^{\mathrm{i}\la_{nm}t}\vp_{nm}^0 \pt_{x_1}\Phi_{n}^{(1)}(1) \Phi_{m}^{(2)}\\
			&= \left(e^{\mathrm{i}\la_{\de,nm}t}-e^{\mathrm{i}\la_{nm}t}\right)\vp_{\de,nm}^0 \pt_{x_1}\Phi_{\de,n}^{(1)}(1) E_\de^0(\Phi_{\de,m}^{(2)}) + e^{\mathrm{i} \la_{nm}t}\left(\vp_{\de,nm}^0-\vp_{nm}^0\right)\pt_{x_1}\Phi_{\de,n}^{(1)}(1) E_\de^0(\Phi_{\de,m}^{(2)})\\
			&\quad + e^{\mathrm{i}\la_{nm}t}\vp_{nm}^0 \left(\pt_{x_1}\Phi_{\de,n}^{(1)}(1)-\pt_{x_1}\Phi_n^{(1)}(1)\right)E_\de^0(\Phi_{\de,m}^{(2)}) + e^{\mathrm{i}\la_{nm}t}\vp_{nm}^0 \pt_{x_1}\Phi_n^{(1)}(1)\left(E_\de^0(\Phi_{\de,m}^{(2)})-\Phi_m^{(2)}\right).
		\end{split}
	\end{equation*}
	By passing to the limit as $\de \to 0^+$, the four terms on the right-hand side vanish uniformly for each $t \in [0,T]$ due to the eigenvalue continuity \eqref{10.20.18}, the convergence of Fourier coefficients \eqref{10.20.19}, the boundary normal derivative trace convergence \eqref{10.20.3}, and the strong eigenfunction convergence guaranteed by Lemma \ref{10.19.L5}, respectively. Since the sum in $D_1$ involves only finitely many terms, there exists a sufficiently small geometric threshold $\de_1 > 0$ such that for all $\de\in (0,\de_1)$, we have 
	\begin{equation*}
		D_1 \leq C\e,
	\end{equation*}
	where the positive constant $C$ remains completely independent of $\de$. Combining the uniform bounds obtained for $D_1, D_2$, and $D_3$, we conclude the proof of the strong convergence of the boundary normal derivatives in $L^2\left(0,T; L^2(\Ga_1^1)\right)$. 
	We have proved the theorem. 
\end{proof}

\begin{lemma}\label{10.21.L1}
	Let $\vp^0\in C_0^\iy(\Om)$ and $\de\in (0,\de_0]$, where $\de_0$ is the geometric threshold parameter defined in \eqref{10.20.6}. Let $\vp_\de$ be the unique weak solution to the perturbed problem \eqref{10.19.6} corresponding to the data $(\vp^0,0)$. Then, the conormal derivative satisfies the uniform admissible boundary trace estimate
	\begin{equation*}
		\left\|\f{\pt (E_\de^0\vp_\de)}{\pt\nu}\right\|_{L^2\left(0,T; L^2(\Ga_1^1\cup \Ga_2^1)\right)}=\left\|\f{\pt (E_\de^0\vp_\de)}{\pt\nu_A}\right\|_{L^2\left(0,T; L^2(\Ga_1^1\cup \Ga_2^1)\right)} \leq C \|\vp^0\|_{H_0^1(\Om;A)},
	\end{equation*}
	where the positive constant $C$ depends exclusively on $\al_1, \al_2$, and $T$, but remains strictly independent of $\de \in (0, \de_0]$.
\end{lemma}

\begin{proof}
	Note that 
	\begin{equation}\label{10.21.1}
		\begin{split}
			\f{\pt (E_\de^0\vp_\de)}{\pt \nu_{A}}
			&=\sum_{n,m=1}^\iy e^{\im \la_{\de,nm}t} \vp_{\de,nm}^0 \left(\pt_{x_1}\Phi_{\de,n}^{(1)}(1)\right) E_\de^0 \Phi_{\de,m}^{(2)} \mbox{ on }\Ga_1^1\ts (0,T),
		\end{split}
	\end{equation}
	it suffices to establish the continuous dependence on the vertical section of the non-degenerate boundary, namely,
	\begin{equation}\label{10.20.18_admit}
		\left\|\f{\pt \vp_\de}{\pt\nu_{A_\de}}\right\|_{L^2\left(0,T; L^2(\Ga_{\de,1}^1)\right)}=\left\|\f{\pt (E_\de^0\vp_\de)}{\pt\nu_A}\right\|_{L^2\left(0,T; L^2(\Ga_1^1)\right)} \leq C \|\vp^0\|_{H_0^1(\Om;A)},
	\end{equation}
	where the positive constant $C$ depends only on $\al_1, \al_2$, and $T$. An analogous argument directly yields the parallel estimate for the horizontal segment $\Ga_2^1$.
	
	Let $\zeta = \zeta(x_1) \in C^\iy(\R),0\leq\zeta\leq 1$ such that 
	\begin{equation*}
		\zeta=0 \mbox{ on } \left(-\iy, \f{1}{4}\right), \quad \zeta=1 \mbox{ on } \left(\f{1}{2},+\iy\right),\quad |\zeta'|+|\zeta''|\leq C \mbox{ on }\R, 
	\end{equation*}
	where the positive constant $C$ is absolute.

	Multiplying $-\zeta\ol{ \pt_{x_1} \vp_\de}-\zeta'\ol\vp_\de$ on the both sides of \eqref{10.19.6} with $f_\de=0$, integration on $Q_\de$, from Lemma \ref{10.19.L1} (i.e., $\vp_\de\in L^2(0,T; H^2(\Om_\de))$ and $\pt_t\vp_\de\in L^2(Q_\de)$), we obtain 
	\begin{equation*}
		\begin{split}
			0
			&=\im \iint_{Q_\de} (\pt_t\vp_\de) \left(-\zeta \ol{\pt_{x_1}\vp_\de} -\zeta'\ol\vp_\de\right)\df x\df t-\iint_{Q_\de} (\Div(A_\de\nabla\vp_\de))\left(-\zeta \ol{\pt_{x_1}\vp_\de} -\zeta'\ol\vp_\de\right)\df x\df t\\
			&=:F_1+F_2.
		\end{split}
	\end{equation*}
	
	Integration by parts on $ [0,T]$, since $\vp_\de=0$ on $\Sigma_\de$, we compute 
	\begin{equation*}
		\begin{split}
			F_1
			&=-\im \left.\int_{\Om_\de} \vp_\de \left(\zeta \ol{\pt _{x_1}\vp_\de}+\zeta'\ol\vp_\de\right)\df x\right|_{t=0}^{t=T}+\im \iint_{Q_\de} \vp_\de \left(\zeta \ol{\pt_{x_1}(\pt_{t}\vp_\de)}+\zeta'\ol{\pt_t\vp_\de}\right)\df x\df t\\
			&=-\im \left.\int_{\Om_\de} \zeta \vp_\de \ol{\pt_{x_1}\vp_\de}\df x\right|_{t=0}^{t=T}-\im\left.\int_{\Om_\de} \zeta'|\vp_\de|^2\df x\right|_{t=0}^{t=T}-\im \iint_{Q_\de} \zeta (\pt_{x_1}\vp_\de) \ol{\pt_t\vp_\de}\df x\df t, 
		\end{split}
	\end{equation*}
	then, from 
	\begin{equation*}
		\begin{split}
			&-2\im \iint_{Q_\de} \zeta (\pt_{x_1}\vp_\de) \ol{\pt_t\vp_\de}\df x\df t+2\im\iint_{Q_\de} \zeta \ol{\pt_{x_1}\vp_\de} \pt_t\vp_\de \df x\df t\\
			&=-\im \iint_{Q_\de} \zeta' (\pt_t\vp_\de) \ol\vp_\de\df x\df t+\im \iint_{Q_\de} \zeta' \ol{\pt_t\vp_\de} \vp_\de\df x\df t\\
			&\hspace{4.5mm}-\im \left.\int_{\Om_\de} \zeta (\pt_{x_1}\vp_\de) \ol\vp_\de\df x\right|_{t=0}^{t=T}+\im\left.\int_{\Om_\de} \zeta \ol{\pt_{x_1}\vp_\de} \vp_\de\df x\right|_{t=0}^{t=T}, 
		\end{split}
	\end{equation*}
	we obtain
	\begin{equation*}
		\begin{split}
			2\Real F_1
			&=\Real\left(- \im \left.\int_{\Om_\de} \zeta \vp_\de \ol{\pt_{x_1}\vp_\de}\df x\right|_{t=0}^{t=T}-\im \iint_{Q_\de} \zeta'(\pt_t\vp_\de)\ol\vp_\de \df x\df t\right). 
		\end{split}
	\end{equation*}
	Note that from $\vp_\de=0$ on $\Sigma_\de$ we obtain 
	\begin{equation*}
		\begin{split}
			-\im \iint_{Q_\de} \zeta'(\pt_t\vp_\de)\ol\vp_\de \df x\df t
			&=-\iint_{Q_\de}\zeta'(\Div(A_\de\nabla\vp_\de))\ol\vp_\de\df x\df t=\iint_{Q_\de}A_\de\nabla\vp_\de \cdot \nabla (\zeta' \ol\vp_\de)\df x\df t,  
		\end{split}
	\end{equation*}
	and then 
	\begin{equation*}
		\begin{split}
			2\Real F_1
			&=\Real\left(- \im \left.\int_{\Om_\de} \zeta \vp_\de \ol{\pt_{x_1}\vp_\de}\df x\right|_{t=0}^{t=T}\right)+\Real  \iint_{Q_\de}A_\de\nabla\vp_\de \cdot \nabla (\zeta'  \ol\vp_\de)\df x\df t. 
		\end{split}
	\end{equation*}
	
	Since $\f{\pt \vp_{\de}}{\pt \nu_{A_\de}}=\nabla \vp_\de\cdot A_\de\nu=[\pt_{x_1}\vp_\de](1,\cdot,\cdot)$ holds on $\Ga_{\de,1}^1\ts (0,T)$, and $\pt_{x_1}\vp_\de=0$ on $(\Ga_{\de,2}^\de\cup \Ga_{\de,2}^1)\ts (0,T)$ by \eqref{10.21.1}, we have 
	\begin{equation*}
		\begin{split}
			F_2
			&=\iint_{\Sigma_\de} (\nabla\vp_\de\cdot A_\de\nu)(\zeta\ol{\pt_{x_1}\vp_\de})\df S\df t-\iint_{Q_\de}A_\de\nabla\vp_\de\cdot \nabla \left(\zeta \pt_{x_1}\ol\vp_\de+\zeta'\ol\vp_\de\right)\df x\df t\\
			&=\iint_{\Ga_{\de,1}^1\ts (0,T)} |\pt_{x_1}\vp_\de|^2\df S\df t-\iint_{Q_\de}\zeta' x_1^{\al_1}|\pt_{x_1}\vp_\de|^2\df x\df t-\iint_{Q_\de} \zeta' \nabla \vp_\de \cdot A_\de\nabla\ol\vp_\de\df x\df t\\
			&\hspace{4.5mm}-\iint_{Q_\de} \zeta'' x_1^{\al_1} (\pt_{x_1}\vp_\de)\ol\vp_\de\df x\df t-\iint_{Q_\de} \zeta \nabla\vp_\de\cdot A_\de \nabla \ol{\pt_{x_1}\vp_\de}\df x\df t. 
		\end{split}
	\end{equation*}
	Together with this and 
	\begin{equation*}
		\begin{split}
			&-\iint_{Q_\de} \zeta \nabla\vp_\de\cdot A_\de \nabla \ol{\pt_{x_1}\vp_\de}\df x\df t-\iint_{Q_\de} \zeta \nabla\ol\vp_\de\cdot A_\de \nabla \pt_{x_1}\vp_\de \df x\df t\\
			&=-\sum_{i=1}^2\iint_{Q_\de} \zeta A_{ii}\f{\pt }{\pt x_1} \left|\f{\pt \vp_\de}{\pt x_i}\right|^2\df x\df t\\
			&=-\iint_{\Ga_{\de,1}^1\ts (0,T)}|\pt_{x_1}\vp_\de|^2\df S\df t+\sum_{i=1}^2\iint_{Q_\de}\left|\f{\pt \vp_\de}{\pt x_i}\right|^2\pt_{x_1}(\zeta A_{ii})\df x\df t
		\end{split}
	\end{equation*}
	by $A=(A_{ij})_{i,j=1,2}$ and $\pt_{x_2}\vp_\de=0$ on $(\Ga_{\de,1}^\de\cup \Ga_{\de,1}^1)\ts (0,T)$, we obtain 
	\begin{equation*}
		\begin{split}
			2\Real F_2
			&=\iint_{\Ga_{\de,1}^1\ts (0,T)}|\pt_{x_1}\vp_\de|^2\df S\df t-2\iint_{Q_\de}\zeta' x_1^{\al_1}|\pt_{x_1}\vp_\de|^2\df x\df t-2\iint_{Q_\de} \zeta' \nabla \vp_\de \cdot A_\de\nabla\ol\vp_\de\df x\df t\\
			&\hspace{4.5mm}-2\Real \iint_{Q_\de} \zeta'' x_1^{\al_1} (\pt_{x_1}\vp_\de)\ol\vp_\de\df x\df t+\sum_{i=1}^2\iint_{Q_\de}\left|\f{\pt \vp_\de}{\pt x_i}\right|^2\pt_{x_1}(\zeta A_{ii})\df x\df t. 
		\end{split}
	\end{equation*} 
	
	Combining  $\Real F_1,\Real F_2$ and Lemma \ref{10.20.L3}, we deduce that 
	\begin{equation*}
		\begin{split}
			\iint_{\Ga_{\de,1}^1\ts (0,T)}|\pt_{x_1}\vp_\de|^2\df S\df t
			&\leq C\left.\int_{\Om_\de} |\vp_\de(t)|^2\df x\right|_{t=0,T}  +C\left.\int_{\Om_\de} |\pt_{x_1}\vp_\de(t)|^2\df x\right|_{t=0,T}\\
			&\hspace{4.5mm}+C\iint_{Q_\de} |\vp_\de|^2\df x\df t+C\iint_{Q_\de}\nabla\vp_\de \cdot A_\de\nabla\vp_\de\df x\df t\\
			&\leq C\|\vp^0\|_{H_0^1(\Om_\de;A_\de)}^2= C\|\vp^0\|_{H_0^1(\Om;A)}^2, 
		\end{split}
	\end{equation*}
	where the positive constants $C$ depend only on $\al_1,\al_2$ and $T$. We prove the lemma. 
\end{proof}

\begin{corollary}\label{10.21.C1}
	Let $\vp^0\in C_0^\iy(\Om)$. Let $\vp$ be the unique weak solution to the perturbed problem \eqref{10.18.11} corresponding to the data $(\vp^0,0)$. Then, 
	\begin{equation}\label{10.21.2}
		\left\|\f{\pt \vp}{\pt\nu}\right\|_{L^2\left(0,T; L^2(\Ga_1^1\cup \Ga_2^1)\right)}=\left\|\f{\pt \vp}{\pt\nu_A}\right\|_{L^2\left(0,T; L^2(\Ga_1^1\cup \Ga_2^1)\right)} \leq C \|\vp^0\|_{H_0^1(\Om;A)},
	\end{equation}
	where the positive constant $C$ depends exclusively on $\al_1, \al_2$, and $T$.
\end{corollary}

\begin{proof}
	From \eqref{10.20.14}, we obtain 
	\begin{equation*}
		\f{\pt (E_\de^0\vp_\de)}{\pt \nu_{A}}\ra \f{\pt\vp}{\pt \nu_A} \mbox{ strongly in } L^2\left(0,T; L^2(\Ga_1^1\cup \Ga_2^1)\right). 
	\end{equation*}
	Combining this and Lemmas \ref{10.20.L2} and \ref{10.21.L1}, we deduce \eqref{10.21.2}. We prove the corollary. 
\end{proof}

\begin{remark}\label{10.21.R_hidden}
	The boundary trace estimate established in Corollary \ref{10.21.C1} cannot be obtained directly via the classical hidden regularity arguments. Indeed, due to the critical boundary degeneracy of the operator $\mcA$, its domain $D(\mcA)$ fails to provide the global $H^2(\Om)$-regularity that is fundamentally required to justify the classical Rellich multiplier method. Consequently, the variational integration-by-parts framework leading to the boundary trace identities cannot be rigorously performed on the weak solutions of the original degenerate problem. Instead, this difficulty is circumvented by introducing a parameterized family of uniformly elliptic domain approximations, performing the multiplier arguments on the non-degenerate subdomains, and subsequently passing to the limit through the strong convergence results established in Section \ref{S3.4}.
\end{remark}

\section{Observability}\label{S4}

In this section, we concentrate on establishing the boundary observability inequalities, focusing sequentially on the parameterized uniformly elliptic equations and the original degenerate Schr\"odinger system.

\subsection{Observability for uniformly approximations}\label{S4.1}

Let $\vp^0\in C_0^\iy(\Om)$, and let $\de\in (0,\de_0]$ with $\de_0$ being the geometric threshold parameter defined in \eqref{10.20.6}. We denote by $\vp_\de$ the unique weak solution to the perturbed homogeneous problem \eqref{10.19.6} corresponding to the initial datum $\vp^0$.

By virtue of the higher-order spectral regularity established in Lemma \ref{10.19.L1}, the state solution exhibits the strong regularities $\vp_\de \in L^2(0,T; H^2(\Om_\de))$ and $\pt_t\vp_\de \in L^2(Q_\de)$. Consequently, the subsequent space-time variational formulations and the sophisticated integration-by-parts procedures are rigorously justified on the parameterized subdomains.

Let the multiplier be 
\begin{equation*}
	x\cdot \nabla\ol \vp_\de +\ol\vp_\de. 
\end{equation*}
Denote 
\begin{equation*}
	\Ga_\de^+=\left\{x\in \pt\Om_\de\colon x\cdot\nu>0\right\},\quad \Sigma_\de^+=\Ga_\de^+\ts (0,T), 
\end{equation*}
where $\nu$ is the outer normal vector of $\pt\Om_\de$. Since
\begin{equation}\label{10.21.3}
	x\cdot\nu=-\de \mbox{ on }\Ga_{\de,1}^\de\cup \Ga_{\de,2}^\de, \quad x\cdot \nu=1 \mbox{ on }\Ga_{\de,1}^1\cup \Ga_{\de,2}^1, 
\end{equation}
we obtain $\Ga_\de^+=\Ga_{\de,1}^1\cup \Ga_{\de,2}^1$. 

Multiplier \eqref{10.19.6} by $x\cdot\nabla\ol\vp_\de +\ol\vp_\de$, integration on $Q_\de$, we obtain 
\begin{equation*}
	\begin{split}
		0=\im \iint_{Q_\de}(\pt_t\vp_\de)(x\cdot \nabla \ol\vp_\de+\ol\vp_\de)\df x\df t-\iint_{Q_\de}[\Div(A_\de\nabla\vp_\de)](x\cdot\nabla\ol\vp_\de +\ol\vp_\de)\df x\df t=:G_1+G_2.
	\end{split}
\end{equation*}

Since $\vp_\de=0$ on $\Sigma_\de$, we obtain 
\begin{equation*}
	\begin{split}
		&2\im\iint_{Q_\de}(\pt_t\vp_\de)(x\cdot \nabla\ol\vp_\de)\df x\df t-2\im \iint_{Q_\de} \ol{\pt_t\vp_\de}(x\cdot \nabla\vp_\de)\df x\df t\\
		&=\im \left.\int_{\Om_\de} \vp_\de(x\cdot \nabla\ol\vp_\de)\df x\right|_{t=0}^{t=T}-\im \left.\int_{\Om_\de} \ol\vp_\de(x\cdot \nabla\vp_\de)\df x\right|_{t=0}^{t=T}\\
		&\hspace{4.5mm}+2\im \iint_{Q_\de} \vp_\de \pt_t\ol\vp_{\de}\df x\df t-2\im \iint_{Q_\de} (\pt_t\vp_\de)\ol\vp_\de\df x\df t,
	\end{split}
\end{equation*}
and then 
\begin{equation*}
	\begin{split}
		2\Real G_1
		&=-\Image  \left.\int_{\Om_\de} \vp_\de(x\cdot \nabla\ol\vp_\de)\df x\right|_{t=0}^{t=T}. 
	\end{split}
\end{equation*}

From $\vp_\de=0$ on $\Sigma_\de$ we obtain $\nabla\vp_\de=(\nabla\vp_\de\cdot\nu)\nu$ on $\Sigma_\de$, and then  
\begin{equation*}
	\begin{split}
		&\iint_{Q_\de}A_\de\nabla\vp_\de \cdot \left[(D^2\ol\vp_\de)x\right]\df x\df t+\iint_{Q_\de} A_\de\nabla\ol\vp_\de\cdot \left[(D^2\vp_\de)x\right]\df x\df t\\
		&=\sum_{i,j=1}^2\iint_{Q_\de} A_{ii}x_j \f{\pt}{\pt x_j}\left|\f{\pt\vp_\de}{\pt x_i}\right|^2 \df x\df t\\
		&=\iint_{\Sigma_\de}(\nabla\vp_\de\cdot A_\de\nabla\ol\vp_\de)(x\cdot\nu)\df S\df t-\sum_{i,j=1}^2\iint_{Q_\de}\left|\f{\pt\vp_\de}{\pt x_i}\right|^2 \f{\pt }{\pt x_j}(A_{ii}x_j)\df x\df t\\
		&=\iint_{\Sigma_\de}\left|\f{\pt \vp_\de}{\pt\nu}\right|^2(A_\de\nu\cdot \nu)(x\cdot\nu)\df S\df t\\
		&\hspace{4.5mm}-2\iint_{Q_\de} \nabla\vp_\de\cdot A_\de\nabla\ol\vp_\de\df x\df t-\sum_{i=1}^2\al_i\iint_{Q_\de} x_i^{\al_i}|\pt_{x_i}\vp_\de|^2\df x\df t.
	\end{split}
\end{equation*}
Together with this and $\nabla\vp_\de=(\nabla\vp_\de\cdot\nu)\nu$ on $\Sigma_\de$ we obtain 
\begin{equation*}
	\begin{split}
		2\Real G_2
		&=-\iint_{\Sigma_\de}\left|\f{\pt\vp_\de}{\pt\nu}\right|^2(A_\de\nu\cdot\nu)(x\cdot\nu)\df S\df t\\
		&\hspace{4.5mm}+2\iint_{Q_\de} \nabla\vp_\de\cdot A_\de\nabla\vp_\de\df x\df t-\sum_{i=1}^2\al_i\iint_{Q_\de} x_i^{\al_i} \left|\pt_{x_i}\vp_\de\right|^2\df x\df t. 
	\end{split}
\end{equation*}

Combining $\Real G_1, \Real G_2$, and \eqref{10.21.3},   we deduce that 
\begin{equation*}
	\begin{split}
		\iint_{\Sigma_\de^+}\left|\f{\pt\vp_\de}{\pt\nu}\right|^2\df S\df t
		&=\de^{1+\al_1}\iint_{\Ga_{\de,1}^\de\ts (0,T)}\left|\f{\pt\vp_\de}{\pt\nu}\right|^2\df S\df t+\de^{1+\al_2}\iint_{\Ga_{\de,2}^\de\ts (0,T)}\left|\f{\pt\vp_\de}{\pt\nu}\right|^2\df S\df t\\
		&\hspace{4.5mm}+\sum_{i=1}^2(2-\al_i)\iint_{Q_\de} x_i^{\al_i}|\pt_{x_i}\vp_\de|^2\df x\df t-\Image\left. \int_{\Om_\de} \vp_\de (x\cdot \nabla\ol\vp_\de)\df x\right|_{t=0}^{t=T}. 
	\end{split}
\end{equation*}
From  Lemma  \ref{10.18.7} we obtain 
\begin{equation*}
	\begin{split}
		\left|\int_{\Om_\de} \vp_\de (x\cdot\nabla\ol\vp_\de)\df x\right|
		&\leq \ga \int_{\Om_\de}|\vp_\de|^2\df x+\f{1}{4\ga}\int_{\Om_\de} |x\cdot\nabla\vp_\de|^2\df x\\
		&\leq \left(\f{4\ga}{(1-\al_i)^2}+\f{1}{2\ga}\right)\int_{\Om_\de} \nabla \vp_\de\cdot A_\de\nabla\ol \vp_\de\df x. 
	\end{split}
\end{equation*}
Choosing $\ga=\f{\sqrt{2}}{4}(1-\al_i)$ we obtain 
\begin{equation*}
	\left|\int_{\Om_\de} \vp_\de (x\cdot \nabla\ol\vp_\de)\df x\right|\leq \f{2\sqrt{2}}{1-\al_i}\int_{\Om_\de} \nabla\vp_\de \cdot A_\de\nabla\ol\vp_\de\df x, \mbox{ for } i=1,2. 
\end{equation*}
Therefore,
\begin{equation*}
	\begin{split}
		\left(\min\{2-\al_1,2-\al_2\}T-\min\left\{\f{4\sqrt{2}}{1-\al_1}, \f{4\sqrt{2}}{1-\al_2}\right\} \right)\|\vp^0\|_{H_0^1(\Om;A)}^2\leq \iint_{\Sigma_\de^+}\left|\f{\pt \vp_\de}{\pt\nu}\right|^2\df S\df t
	\end{split}
\end{equation*}
by Lemma \ref{10.20.L3}, this shows the following theorem. 

\begin{theorem}\label{10.21.T2}
	Let $\vp^0\in C_0^\iy(\Om)$ and $\de_0>0$ be defined in \eqref{10.20.6}. Denote $\vp_\de\ (0<\de<\de_0)$ be the weak solution of \eqref{10.19.6} with respect to $(\vp^0,0)$. Then
	\begin{equation}\label{10.21.6}
		\begin{split}
			\left(\min\{2-\al_1,2-\al_2\}T-\min\left\{\f{4\sqrt{2}}{1-\al_1}, \f{4\sqrt{2}}{1-\al_2}\right\} \right)\|\vp^0\|_{H_0^1(\Om;A)}^2\leq \iint_{\Sigma_\de^+}\left|\f{\pt \vp_\de}{\pt\nu}\right|^2\df S\df t. 
		\end{split}
	\end{equation}
\end{theorem}

\subsection{Observability for the original equation}\label{S4.2}

Denote $\Sigma^+=\Ga^+\ts (0,T)$ with $\Ga^+=\Ga_1^1\cup \Ga_2^1$. Note that on the non-degenerate active boundary $\Ga_1^1 \cup \Ga_2^1$, since the coordinates satisfy $x_1 = 1$ and $x_2 = 1$ respectively, the matrix profiles collapse to $1^{\al_i} = 1$. Consequently, the conormal derivative naturally coincides with the standard normal derivative, i.e., $\f{\pt \vp}{\pt\nu_A} = \f{\pt \vp}{\pt\nu} = \pt_{x_1}\vp(1, \cdot, \cdot)$ on $\Ga_1^1 \times (0,T)$, and $\f{\pt \vp}{\pt\nu_A} = \f{\pt \vp}{\pt\nu} = \pt_{x_2}\vp(1, \cdot, \cdot)$ on $\Ga_2^1 \times (0,T)$.

\begin{theorem}\label{10.21.T1}
	Let $\al_1, \al_2 \in (0,1)$ with $\al_1 \leq \al_2$. Define the minimal observability time $T_0$ by
	\begin{equation}\label{10.21.4}
		T_0 = \f{4\sqrt{2}}{(1-\al_1)(2-\al_2)}.
	\end{equation}
	Then, for any $T > T_0$, there exists a positive observability constant $C$ depending only on $\al_1, \al_2$, and $T$, such that the unique weak solution $\vp$ to the original degenerate Schr\"odinger equation \eqref{10.18.11} satisfies the reverse observability inequality:
	\begin{equation}\label{10.21.5}
		\|\vp^0\|_{H_0^1(\Om;A)}^2 \leq C \iint_{\Sigma^+}\left|\f{\pt \vp}{\pt\nu}\right|^2\df S\,\df t
	\end{equation}
	for all initial data $\vp^0 \in C_0^\iy(\Om)$.
\end{theorem}

\begin{proof}
	Let $\de_0$ be defined in \eqref{10.20.6}. Let $\de\in (0,\de_0]$, denote $\vp_\de$ be the weak solution of \eqref{10.19.6} with respect to $(\vp^0,0)$. Then we have \eqref{10.21.6}. Note that 
	\begin{equation*}
		\f{4\sqrt{2}}{1-\al_1}\leq \f{4\sqrt{2}}{1-\al_2}, 
	\end{equation*}
	we get 
	\begin{equation*}
		\begin{split}
			\left((2-\al_2)T-\f{4\sqrt{2}}{1-\al_1}\right)\|\vp^0\|_{H_0^1(\Om;A)}^2\leq \iint_{\Sigma_\de^+}\left|\f{\pt\vp_\de}{\pt\nu}\right|^2\df S\df t. 
		\end{split}
	\end{equation*}
	Combining this and \eqref{10.20.14},  we obtain 
	\begin{equation*}
		\begin{split}
			\left((2-\al_2)T-\f{4\sqrt{2}}{1-\al_1}\right)\|\vp^0\|_{H_0^1(\Om;A)}^2\leq \iint_{\Sigma^+}\left|\f{\pt\vp}{\pt\nu}\right|^2\df S\df t 
		\end{split}
	\end{equation*}
	for all $T>T_0$, where the positive constant $T_0$ defined in \eqref{10.21.4}.
	This proves \eqref{10.21.5}. 
\end{proof}

Now, we prove the equation of \eqref{10.18.11} with initial  data $\vp^0\in H_0^1(\Om;A)$. 

Since $C_0^\iy(\Om)$ is dense in $H_0^1(\Om;A)$, there exists a sequence $\{\vp_n^0\}_{n\in\N^*}\s C_0^\iy(\Om)$ such that 
\begin{equation}\label{10.22.1}
	\vp_n^0\ra \vp^0 \mbox{ strongly in } H_0^1(\Om;A). 
\end{equation}
Denote $\vp_n\ (n\in\N^*)$ be the weak solution of \eqref{10.18.11} with respect to $(\vp_n^0,0)$. From Theorem \ref{10.21.T1}, we obtain for all $T>T_0$,  
\begin{equation}\label{10.22.2}
	\|\vp_n^0\|_{H_0^1(\Om;A)}^2\leq C\iint_{\Sigma^+}\left|\f{\pt \vp_n}{\pt \nu}\right|^2\df S\df t, 
\end{equation}
where the positive constant $C$ depends only on $\al_1,\al_2$ and $T$, and $T_0$ is defined in \eqref{10.21.4}.

Next, we pass to the limit as $n \to \iy$ on both sides of \eqref{10.22.2}. On the left-hand side, the strong convergence in \eqref{10.22.1} immediately implies that \begin{equation*}
	\lim_{n\to\iy} \|\vp_n^0\|_{H_0^1(\Om;A)}^2 = \|\vp^0\|_{H_0^1(\Om;A)}^2. 
\end{equation*}
On the right-hand side, by virtue of the linear continuous dependence of the boundary conormal derivative on the initial states (guaranteed by the admissibility inequality in Corollary \ref{10.21.C1}), the mapping 
\begin{equation*}
	\vp^0 \mapsto \f{\pt \vp}{\pt\nu_A}
\end{equation*}
is a bounded linear operator from $H_0^1(\Om;A)$ into $L^2\left(0,T; L^2(\Ga_1^1\cup \Ga_2^1)\right)$. Consequently, we deduce that
\begin{equation*}
	\f{\pt \vp_n}{\pt\nu_A} \ra \f{\pt \vp}{\pt\nu_A} \text{ strongly in } L^2\left(0,T; L^2(\Sigma^+)\right) \text{ as } k \to \iy.
\end{equation*}
This strong trace convergence allows us to safely pass to the limit on the right-hand side of \eqref{10.22.2}. We therefore conclude that the reverse observability inequality
\begin{equation*}
	\|\vp^0\|_{H_0^1(\Om;A)}^2 \leq C \iint_{\Sigma^+}\left|\f{\pt \vp}{\pt\nu_A}\right|^2\df S\,\df t
\end{equation*}
holds for all initial data $\vp^0 \in H_0^1(\Om;A)$, where the positive constant $C$ depends only on $\al_1,\al_2$ and $T$. This proves the following  theorem, it is one of main results of this paper. 

\begin{theorem}\label{10.22.T1}
	Let $\al_1, \al_2 \in (0,1)$ with $\al_1 \leq \al_2$. Define the minimal observability time $T_0$ by
	\begin{equation*} 
		T_0 = \f{4\sqrt{2}}{(1-\al_1)(2-\al_2)}.
	\end{equation*}
	Then, for any $T > T_0$, there exists a positive observability constant $C$ depending only on $\al_1, \al_2$, and $T$, such that the unique weak solution $\vp$ to the original degenerate Schr\"odinger equation \eqref{10.18.11} satisfies the reverse observability inequality:
	\begin{equation*}
		\|\vp^0\|_{H_0^1(\Om;A)}^2 \leq C \iint_{\Sigma^+}\left|\f{\pt \vp}{\pt\nu}\right|^2\df S\,\df t
	\end{equation*}
	for all initial data $\vp^0 \in H_0^1(\Om;A)$.
\end{theorem}

\section{Appendix}\label{S5}

In this section, we will prove Theorems \ref{10.18.T1}  and \ref{10.19.T1}. 

\begin{proof}[Proof of Theorem \ref{10.18.T1}]
	We prove this theorem by the following steps. 
	
	{\it Step 1}. Galerkin method. 
	
	Let $m\in\N^*$. Define 
	\begin{equation}\label{09.16.5}
		\vp^m=\sum_{n=1}^m \vp_n^m(t) \Phi_n, \quad f^m=\sum_{n=1}^m f_n^m(t)\Phi_n \mbox{ with } f_n^m(t)=(f,\Phi_n)_{L^2(\Om)}, 
	\end{equation}
	where $\vp_n^m(t), n=1,\cdots, m$ is the solution of the ordinary differential equation 
	\begin{equation}\label{09.16.6}
		\im \pt_t\vp_n^m(t)+\la_n \vp_n^m(t)=f_n^m(t),   \ n=1,\cdots, m, \quad   \pt_t\vp_n^m\equiv \f{\df\vp_n^m}{\df t},
	\end{equation}
	with initial data
	\begin{equation}\label{09.16.7}
		\vp_n^m(0)=\vp_n^0, \mbox{ with } \vp_n^0=(\vp^0,\Phi_n)_{L^2(\Om)}, \ n=1,\cdots, m. 
	\end{equation}
	Indeed, the equation \eqref{09.16.6} is equivalent to the following equation 
	\begin{equation}\label{09.16.8}
		\im(\pt_t\vp^m, \Phi_n)_{L^2(\Om)}+a(\vp^m, \Phi_n)=(f^m, \Phi_n)_{L^2(\Om)}, \ n=1,\cdots, m, 
	\end{equation}
	where $a(\cdot,\cdot)$ and $\Phi_n$ are defined in \eqref{10.18.12} and \eqref{10.18.15}. 
	
	{\it Step 2}. Energy estimate. 
	
	Multiplying $\ol{\vp_n^m}$ on the both sides of \eqref{09.16.6}, summing $n=1,\cdots, m$, we obtain 
	\begin{equation}\label{09.16.9}
		\begin{split}
			\im \sum_{n=1}^m (\pt_t\vp_n^m)\ol{\vp_n^m}+\sum_{n=1}^m \la_n \vp_n^m \ol{\vp_n^m}=\sum_{n=1}^m f_n^m\ol{\vp_n^m}. 
		\end{split}
	\end{equation}
	Taking conjugate of \eqref{09.16.9}, we obtain 
	\begin{equation*}
		-\im \sum_{n=1}^m \ol{(\pt_t\vp_n^m)}{\vp_n^m}+\sum_{n=1}^m \la_n \vp_n^m \ol{\vp_n^m}=\sum_{n=1}^m \ol{f_n^m}{\vp_n^m}.
	\end{equation*}
	Then 
	\begin{equation*}
		\begin{split} 
			\f{\df}{\df t}\|\vp^m(t)\|_{L^2(\Om)}^2
			&=\sum_{n=1}^m \pt_t\left(\vp_n^m \ol{\vp_n^m}\right)=2\Image \sum_{n=1}^m f_n^m \ol{\vp_n^m}\\
			&=2\Image (f^m, \vp^m)_{L^2(\Om)}\leq \|f^m\|_{L^2(\Om)}^2+\|\vp^m\|_{L^2(\Om)}^2. 
		\end{split} 
	\end{equation*}
	The Gronwall's inequality \cite[Chapter B.2 (j), p.~708]{Evans}  implies 
	\begin{equation*}
		\begin{split}
			\|\vp^m(t)\|_{L^2(\Om)}^2\leq e^t \left(\|\vp^m(0)\|_{L^2(\Om)}^2+\int_0^t \|f^m(s)\|_{L^2(\Om)}^2\df s\right), \mbox{ for all } t\in [0,T]. 
		\end{split}
	\end{equation*}
	Hence, for all $m\in\N^*$, 
	\begin{equation}\label{09.16.10}
		\esssup_{t\in [0,T]}\|\vp^m(t)\|_{L^2(\Om)}^2\leq e^T\left(\|\vp^0\|_{L^2(\Om)}^2+\|f\|_{L^2(Q)}^2\right). 
	\end{equation}
	
	Multiplying \eqref{09.16.6} by $\ol{\pt_t\vp_n^m}$, summing $n=1,\cdots, m$, we get
	\begin{equation}\label{09.16.12}
		\begin{split}
			\im \sum_{n=1}^m |\pt_t\vp_n^m|^2+\sum_{n=1}^m \la_n \vp_n^m\ol{\pt_t\vp_n^m}=\sum_{n=1}^m f_n^m \ol{\pt_t \vp_n^m}. 
		\end{split}
	\end{equation}
	Taking conjugate of \eqref{09.16.12}, we obtain 
	\begin{equation*}
		\begin{split}
			-\im \sum_{n=1}^m |\pt_t\vp_n^m|^2+\sum_{n=1}^m \la_n \ol{\vp_n^m}\pt_t\vp_n^m=\sum_{n=1}^m \ol{f_n^m}\pt_t\vp_n^m. 
		\end{split}
	\end{equation*}
	Then 
	\begin{equation}\label{09.16.13}
		\begin{split}
			\f{\df}{\df t}a(\vp^m,\vp^m)
			&=2\Real \sum_{n=1}^m \la_n \vp_n^m \ol{\pt_t\vp_n^m}=2\Real \sum_{n=1}^m f_n^m\ol{\pt_t\vp_n^m}\\
			&=2\Real (f^m,\pt_t\vp^m)_{L^2(\Om)}\leq \|f^m\|_{H_0^1(\Om;w)}^2+\|\pt_t\vp^m\|_{H^{-1}(\Om;w)}^2. 
		\end{split}
	\end{equation}
	Note that from \eqref{10.18.7} and  \eqref{09.16.5} and $a(\Phi_n,v)=\la_n(\Phi_n,v)_{L^2(\Om)}$ for all $v\in H_0^1(\Om;w)$, we obtain 
	\begin{equation}\label{09.16.16}
		\begin{aligned}
			&\|\partial_t\varphi^m\|_{H^{-1}(\Omega;w)}\\
			&=
			\sup_{\|v\|_{H_0^1(\Omega;w)}\leq 1}
			\left|
			\langle \partial_t\varphi^m,v\rangle_{H^{-1}(\Omega;w),H_0^1(\Omega;w)}
			\right| =
			\sup_{\|v\|_{H_0^1(\Omega;w)}\leq 1}
			\left|
			(\partial_t\varphi^m,v)_{L^2(\Omega)}
			\right|
			\\
			&=
			\sup_{\|v\|_{H_0^1(\Omega;w)}\leq 1}
			\left|
			\sum_{n=1}^{m}
			\partial_t\varphi_n^m
			(\Phi_n,v)_{L^2(\Omega)}
			\right| =
			\sup_{\|v\|_{H_0^1(\Omega;w)}\leq 1}
			\left|
			\sum_{n=1}^{m}
			\lambda_n^{-1}
			\partial_t\varphi_n^m
			a(\Phi_n,v)
			\right|
			\\
			&\leq
			\sup_{\|v\|_{H_0^1(\Omega;w)}\leq 1}
			\left(
			\sum_{n=1}^{m}
			\lambda_n^{-1}
			|\partial_t\varphi_n^m|^2
			\right)^{\frac12}
			\left(
			\sum_{n=1}^{m}
			\lambda_n^{-1}
			|a(\Phi_n,v)|^2
			\right)^{\frac12} \leq
			\left(
			\sum_{n=1}^{m}
			\lambda_n^{-1}
			|\partial_t\varphi_n^m|^2
			\right)^{\frac12}.
		\end{aligned}
	\end{equation}
	Combining this and \eqref{09.16.13} and \eqref{09.16.6}, we get
	\begin{equation*}
		\begin{split}
			\f{\df}{\df t}a(\vp^m,\vp^m)
			&\leq \|f^m\|_{H_0^1(\Om;w)}^2+\sum_{n=1}^m |\pt_t\vp_n^m|^2\la_n^{-1}\\
			&\leq \|f^m\|_{H_0^1(\Om;w)}^2+2\sum_{n=1}^m \la_n|\vp_n^m|^2+2\sum_{n=1}^m \la_1^{-1} |f_n^m|^2 \leq C\|f^m\|_{H_0^1(\Om;w)}^2+2a(\vp^m,\vp^m), 
		\end{split}
	\end{equation*}
	where the positive constant $C$ depends only on $\al_1$ and $\al_2$. 
	By the Gronwall's inequality, we obtain 
	\begin{equation*}
		a(\vp^m(t),\vp^m(t))\leq e^{2t}\left( \|\vp^m(0)\|_{H_0^1(\Om;w)}^2+C\int_0^t \|f^m(s)\|_{H_0^1(\Om;w)}^2\df s\right), \mbox{ for all } t\in [0,T], 
	\end{equation*}
	and then 
	\begin{equation}\label{09.16.15}
		\begin{split}
			\esssup_{t\in [0,T]}\|\vp^m(t)\|_{H_0^1(\Om;w)}^2\leq e^{2T}\left(\|\vp^0\|_{H_0^1(\Om;w)}^2+C\|f\|_{L^2(0,T; H_0^1(\Om;w))}^2\right), 
		\end{split}
	\end{equation}
	where the positive constant $C$ depends only on $\al_1$ and $\al_2$. Combining this and \eqref{09.16.6} and  \eqref{09.16.16}, we obtain 
	\begin{equation}\label{09.16.17}
		\begin{split}
			\|\pt_t\vp^m\|_{L^2(0,T; H^{-1}(\Om;w))}^2
			&\leq 2\int_0^Ta(\vp^m,\vp^m)\df t+2\la_1^{-1}\int_0^T \|f^m\|_{L^2(\Om)}^2\df t\\
			&\leq Ce^{2T}(1+T)\left(\|\vp^0\|_{H_0^1(\Om;w)}^2+\|f\|_{L^2(Q)}^2\right), 
		\end{split}
	\end{equation}
	where the positive constant $C$ depends only on $\al_1$ and $\al_2$. 
	
	{\it Step 3}. Approximation or convergence. 
	
	Combining \eqref{09.16.10} and \eqref{09.16.15} and \eqref{09.16.17}, there exists a subsequence of $\{\vp^m\}_{m\in\N^*}$, still denoted by itself, and 
	\begin{equation}\label{09.16.18}
		\wt\vp\in L^\iy(0,T; H_0^1(\Om;w)), \mbox{ with } \pt_t\wt\vp\in L^2(0,T; H^{-1}(\Om;w)), 
	\end{equation} 
	such that 
	\begin{equation}\label{09.16.19}
		\begin{split}
			\vp^m
			&\ra \wt\vp \mbox{ weak-* in } L^\iy(0,T; H_0^1(\Om;w)), \\
			\pt_t\vp^m
			&\ra \pt_t\wt\vp \mbox{ weakly in } L^2(0,T; H^{-1}(\Om;w)). 
		\end{split}
	\end{equation}
	Moreover, 
	\begin{equation}\label{09.16.20}
		\begin{split}
			&\esssup_{t\in [0,T]}\|\wt\vp(t)\|_{H_0^1(\Om;w)}^2+\|\pt_t\wt\vp\|_{L^2(0,T; H^{-1}(\Om;w))}^2\\
			&\leq Ce^{2T}(1+T)\left(\|\vp^0\|_{H_0^1(\Om;w)}^2+\|f\|_{L^2(0,T; H_0^1(\Om;w))}^2\right), 
		\end{split}
	\end{equation}
	where the positive constant $C$ depends only on $\al_1$ and $\al_2$.  
	Moreover, 
	\begin{equation}\label{09.16.24}
		\wt\vp\in C([0,T]; L^2(\Om)). 
	\end{equation}
	
	{\it Step 4}. We prove $\wt\vp$ is a weak solution of \eqref{10.18.24}. 
	
	Let $k\in\N^*$. Denote 
	\begin{equation}\label{09.16.21}
		\psi^k(x,t)=\sum_{n=1}^k \psi_n^k(t) \Phi_n(x), \mbox{ where $\psi_n^k(t)$ are smooth functions for }  n=1,\cdots, k. 
	\end{equation}
	Multiply $\psi_n^k(t)$ on the both sides of \eqref{09.16.8}, integrating on $[0,T]$, for all $m\geq k$, we obtain 
	\begin{equation}\label{09.16.23}
		\begin{split}
			\im \int_0^T (\pt_t\vp^m, \psi^k)_{H^{-1}(\Om;w), H_0^1(\Om;w)}\df t+\int_0^T a(\vp^m, \psi^k)\df t=\int_0^T (f^m, \psi^k)_{L^2(\Om)}\df t. 
		\end{split}
	\end{equation}
	Together with this and \eqref{09.16.19}, we get
	\begin{equation*}
		\begin{split}
			\im \int_0^T \lg\pt_t\wt\vp, \psi^k\rg_{H^{-1}(\Om;w), H_0^1(\Om;w)}\df t+\int_0^T a(\wt\vp, \psi^k)\df t=\int_0^T (f, \psi^k)_{L^2(\Om)}\df t. 
		\end{split}
	\end{equation*}
	Since $
	C^\infty([0,T])
	\otimes
	\mathrm{span}\{\Phi_n\colon n\in\mathbb N^\ast\}$ 
	is dense in
	$L^2(0,T;H_0^1(\Om;w))$,
	we obtain
	\begin{equation}\label{09.16.22}
		\im \int_0^T \lg\pt_t\wt\vp, \psi\rg_{H^{-1}(\Om;w), H_0^1(\Om;w)}\df t+\int_0^T a(\wt\vp, \psi)\df t=\int_0^T (f, \psi)_{L^2(\Om)}\df t 
	\end{equation}
	for all $\psi\in L^2(0,T; H_0^1(\Om;w))$. 
	
	For each $\zeta\in C_0^\iy(0,T)$ and $v\in H_0^1(\Om;w)$, from \eqref{09.16.22}, we obtain 
	\begin{equation*} 
		\im \int_0^T\zeta  \lg\pt_t\wt\vp, v\rg_{H^{-1}(\Om;w), H_0^1(\Om;w)}\df t+\int_0^T\zeta  a(\wt\vp, v)\df t=\int_0^T \zeta (f, v)_{L^2(\Om)}\df t. 
	\end{equation*}
	This shows 
	\begin{equation*}
		\lg \im \pt_t\wt\vp, v\rg_{H^{-1}(\Om;w), H_0^1(\Om;w)}+a(\wt\vp, v)=(f, v)_{L^2(\Om)}
	\end{equation*}
	for a.e.~$t\in [0,T]$ and all $v\in H_0^1(\Om;w)$. 
	
	Now, we show $\wt\vp(0)=\vp^0$. 
	
	From \eqref{09.16.23}, we obtain 
	\begin{equation*}
		-\im\int_0^T(\vp^m, \pt_t\psi^k)_{L^2(\Om)}\df t+\int_0^Ta(\vp^m, \psi^k)\df t=\int_0^T(f^m,\psi^k)_{L^2(\Om)}\df t+(\im\vp^m(0),\psi^k(0))_{L^2(\Om)}
	\end{equation*}
	for all $\psi^k$ in \eqref{09.16.21} satisfying $\psi^k(T)=0$. Together with this and \eqref{09.16.7} and \eqref{09.16.19}, we get
	\begin{equation*}
		\begin{split}
			-\im\int_0^T(\wt \vp, \pt_t\psi^k)_{L^2(\Om)}\df t+\int_0^Ta(\wt \vp, \psi^k)\df t=\int_0^T(f,\psi^k)_{L^2(\Om)}\df t+(\im\vp^0,\psi^k(0))_{L^2(\Om)}
		\end{split}
	\end{equation*}
	for all $\psi^k$ in \eqref{09.16.21} satisfying $\psi^k(T)=0$. Moreover, 
	\begin{equation}\label{09.16.25}
		\begin{split}
			-\im\int_0^T(\wt \vp, \pt_t\psi)_{L^2(\Om)}\df t+\int_0^Ta(\wt \vp, \psi)\df t=\int_0^T(f,\psi)_{L^2(\Om)}\df t+(\im\vp^0,\psi(0))_{L^2(\Om)}
		\end{split}
	\end{equation}
	for all $\psi\in C^1([0,T]; H_0^1(\Om;w))$ in \eqref{09.16.21} satisfying $\psi(T)=0$. From \eqref{09.16.24} and  \eqref{09.16.22}, we obtain 
	\begin{equation*}
		-\im \int_0^T( \wt\vp, \pt_t\psi)_{L^2(\Om)}\df t+\int_0^T a(\wt\vp, \psi)\df t=\int_0^T (f,\psi)_{L^2(\Om)}\df t+(\im \wt\vp(0), \psi(0))_{L^2(\Om)}
	\end{equation*}
	for all $\psi\in C^1([0,T]; H_0^1(\Om;w))$ in \eqref{09.16.21} satisfying $\psi(T)=0$. Together with this and \eqref{09.16.25}, we deduce that 
	\begin{equation*}
		\wt\vp(0)=\vp^0. 
	\end{equation*}
	
	{\it Step 5}. Uniqueness. 
	
	Since the equation \eqref{10.18.24} is linear, we only need to show: if $\vp^0=f=0$, then $\vp=0$. Taking $v=\vp$ as a test function in the weak formulation, we obtain 
	\begin{equation}\label{09.16.26}
		\im \int_0^T\lg \pt_t\vp,  \vp\rg_{H^{-1}(\Om;w), H_0^1(\Om;w)}\df t+\int_0^Ta(\vp, \vp)\df t=0. 
	\end{equation}
	Taking conjugate of \eqref{09.16.26}, we get
	\begin{equation*}
		-\im \int_0^T \lg\ol{\pt_t\vp},\ol \vp\rg_{H^{-1}(\Om;w), H_0^1(\Om;w)}\df t+\int_0^T a(\ol\vp,\ol  \vp)\df t=0. 
	\end{equation*}
	Combining this and \eqref{09.16.26} and $a(\vp,\vp)=\int_\Om \nabla \vp \cdot A\nabla\ol\vp\df x=a(\ol\vp,\ol\vp)$, we deduce that 
	\begin{equation*}
		\int_0^t\f{\df}{\df t}\|\vp(t)\|_{L^2(\Om)}^2\df t=\|\vp(t)\|_{L^2(\Om)}^2, \mbox{ for all }t\in [0,T]. 
	\end{equation*}
	Hence $\vp=0$. 
	
	{\it Step 6}. From Steps 4 and 5, we get the existence and uniqueness of weak solution for \eqref{10.18.24}. From \eqref{09.16.20} we obtain \eqref{10.19.1}. We complete the proof of this theorem. 
	\end{proof}
	
	\begin{proof}[Proof of Theorem \ref{10.19.T1}]
	The proof is similar to that of Theorem \ref{10.18.T1} by Remarks \ref{10.19.R3} and \ref{10.19.R4}. 
	Indeed, the Galerkin approximation, the energy estimates, and the
	compactness argument can be performed in the same way on
	$\Omega_\delta$. Therefore, we only need to verify that the estimates
	are uniform with respect to $\delta$.
	
	Since $\Omega_\delta\s  \Omega$ and $\delta\in(0,\f{1}{4}]$, the estimates 
	used in the proof of Theorem \ref{10.18.T1}, including the Gronwall
	inequality, the estimate of $
	\|\partial_t\vp_\delta^m\|_{H^{-1}(\Omega_\delta;w_\delta)}$, 
	and the corresponding spectral estimate, remain valid on
	$\Omega_\delta$. Hence, passing to the limit as in the proof of
	Theorem \ref{10.18.T1}, we obtain a unique weak solution satisfying
	\begin{equation*} 
		\begin{split}
			&\esssup_{t\in[0,T]}
			\|\vp_\delta(t)\|_{H_0^1(\Omega_\delta;w_\delta)}
			+\|\partial_t\vp_\delta\|_
			{L^2(0,T;H^{-1}(\Omega_\delta;w_\delta))}
			\\
			&\leq C\left(
			\|\vp_\delta^0\|_{H_0^1(\Omega_\delta;w_\delta)}
			+\|f_\delta\|_{L^2(0,T;H_0^1(\Omega_\delta;w_\delta))}
			\right),
		\end{split}
	\end{equation*}
	where $C$ depends only on $\alpha_1,\al_2$ and $T$, but is independent of
	$\delta$.
	\end{proof}

\end{document}